\documentclass[12pt,a4paper,usenames,dvipsnames]{article}

\usepackage{fullpage}
\usepackage{psfrag}
\usepackage{amsmath}
\usepackage{fancybox}
\usepackage{amsfonts}
\usepackage{amsthm}
\usepackage{amssymb}
\usepackage{comment}
\usepackage{graphicx}

\usepackage[normalem]{ulem} 
\usepackage{tikz}
\usetikzlibrary{shapes.geometric}
\usepackage{enumerate}
\usepackage{MnSymbol}
\usepackage{bm}

\DeclareMathOperator*{\argmin}{argmin}

\usepackage{algorithm}
\usepackage{setspace}
\usepackage{algpseudocode}

\usepackage{booktabs}
\usepackage{blkarray}
\usepackage{longtable}
\usepackage{colortbl}
\usepackage{multirow}

\usepackage{pifont}
 \usepackage{xargs} 

 \theoremstyle{plain}
    \newtheorem{theorem}{Theorem}[section]
    \newtheorem{corollary}[theorem]{Corollary}
    \newtheorem{lemma}[theorem]{Lemma}
    \newtheorem{proposition}[theorem]{Proposition}
    
 \theoremstyle{definition}
    
    \newtheorem{remark}[theorem]{Remark}
 \theoremstyle{remark}
    \newtheorem{example}[theorem]{Example}

 \usepackage[english]{babel}

\usepackage[%
   bookmarks=true,
   pagebackref=false, 
   colorlinks,
   linkcolor=blue,
   anchorcolor=blue,
   citecolor=blue,
   filecolor=blue,
   urlcolor=blue
   ]{hyperref}

\usepackage{bm}
\usepackage{xspace}

\usepackage{color}

\DeclareMathOperator{\supp}{supp}

\DeclareMathOperator{\Span}{span}

\newcommand{\loc}{\mathrm{loc}}
\newcommand{\fineonly}{\mathrm{fine\text{-}only}} 

\newcommand{\co}{\textrm{co}}

\newcommand\MM{\mathcal M}
\newcommand\Tau{\mathcal T}

\newcommand\QQ{\mathcal Q}
\newcommand\GG{\mathcal G}
\newcommand\HH{\mathcal H}

\newcommand{\RR}{\mathbb R}

\renewcommand{\Xi}{{\bm\xi}}

\renewcommand{\nu}{{\bm\nu}}

\usepackage{authblk}
\newcommand{\norm}[1]{{\left\vert\kern-0.25ex\left\vert\kern-0.25ex\left\vert #1 
    \right\vert\kern-0.25ex\right\vert\kern-0.25ex\right\vert}}

\newcommand{\children}[1]{{\bf children}({#1})}
\newcommand{\parents}[1]{{\bf parents}({#1})}
 
\begin{document}

\title{\bf Algebraic coarsening of hierarchical splines with guaranteed error control}
\author{Marcelo Actis, Silvano Figueroa \& Eduardo M. Garau\thanks{Corresponding author: Eduardo M. Garau (egarau@santafe-conicet.gov.ar)}}
\affil{\small Universidad Nacional del Litoral, \\
    Consejo Nacional de Investigaciones Cient\'ificas y T\'ecnicas, \\
    FIQ, Santa Fe, Argentina}

\maketitle

\begin{abstract}
    We present a purely algebraic coarsening operator for hierarchical spline
    spaces, obtained by inverting locally the refinement matrix that relates the
    coarsened basis to the current one. Only the B-splines that interact with the coarsened region
    are involved, so that a coarsening step reduces to one small least-squares
    problem followed by explicit corrections of the surviving coefficients,
    without mass matrices or numerical quadrature. The operator is a left
    inverse of the refinement matrix, hence exact whenever the fine spline belongs
    to the coarse space, and leaves the spline unchanged away from the coarsened
    region. The residual of the least-squares problem bounds the $L^\infty$-error with no unknown constants and
    is available before the mesh is modified, which yields an adaptive multi-level
    algorithm that guarantees a prescribed tolerance. The construction carries over to
    truncated hierarchical B-splines, where the corrections disappear.
    In numerical experiments, the bound is sharp within an order of magnitude, the local transfer is nearly as accurate as a global $L^2$-projection, and much of the refinement left behind a moving front is removed without degrading the solution.
\end{abstract}

\begin{quote}\small
\end{quote}



\section{Introduction}
\label{sec:introduction}

Isogeometric analysis~\cite{HCB2005} uses the same spline bases, typically tensor-product B-splines or NURBS,  to describe the geometry and
to approximate the solution of a partial differential equation.
Adaptivity, however, requires spline spaces that can be refined locally,
and this is exactly what the tensor-product construction does not offer: a knot
inserted to resolve a localized feature propagates across the whole
parametric domain.
Hierarchical B-splines (HB-splines), introduced by Kraft~\cite{Kraft} and
brought into the isogeometric setting in~\cite{VGJS2011}, overcome this
limitation by superposing nested tensor-product spaces over a nested sequence of
subdomains, so that fine functions are only introduced where they are needed.
Truncated hierarchical B-splines (THB-splines)~\cite{GJS12} recover the
partition of unity, lead to sparser and, in general, better conditioned
system matrices~\cite{GJKMSS2016}, and are widely used in adaptive isogeometric
methods.
 
Adaptivity involves two operations on the mesh, and in the hierarchical setting
they are far from symmetric. Refinement is exact: since the spaces are nested,
the discrete solution on the coarse mesh also belongs to the finer space,
and its new coefficients follow from a sparse matrix--vector product.
Coarsening is not: the coarse space is a proper subspace of the fine one, and a
fine spline has in general no exact representation in it, so that the solution
must be \emph{projected} and some information is lost. For
stationary problems this asymmetry is harmless, since refinement alone yields
optimal convergence rates and the standard adaptive loop
\textsf{solve}\,$\to$\,\textsf{estimate}\,$\to$\,\textsf{mark}\,$\to$\,%
\textsf{refine} has no coarsening step at all~\cite{BGGPV2022}. Coarsening becomes
essential for evolution problems and moving-feature simulations, where the
region that requires resolution travels through the domain and the mesh must be
un-refined behind it if the number of degrees of freedom is to remain bounded. It has a long tradition in adaptive finite element methods,
where it was first used to prove optimal complexity~\cite{BDD2004}, and it is
indispensable for parabolic problems, both with finite elements~\cite{KMSS2012} and in isogeometric analysis~\cite{GGMV2019}. For hierarchical splines, the mesh
side of coarsening---which cells can be un-refined, and how the hierarchy and
the basis are updated---is settled in~\cite{GV18}, where refinement and
coarsening are defined as inverse operations of each other. The question
addressed in this paper is the other half: how to transfer the discrete
solution to a coarsened space, how to quantify what is lost in doing so, and how
to use that information to decide where, and how far, the mesh may be coarsened.
 
For the transfer itself, the most common choice is a global $L^2$-projection onto the coarsened
space~\cite{CGRV2019,BGRTV23}. It is robust, optimal in the $L^2$-norm and easy to combine with mesh
grading, but it requires a mass matrix assembled by numerical
quadrature and a global linear system solved at every coarsening step. Moreover, since the inverse of the mass matrix couples all the degrees of freedom, a
strictly local un-refinement perturbs the solution also where the mesh does not change. The projection also gives no indication of the accuracy lost when the fine solution does not belong to the coarse space: Carraturo et al., who for efficiency perform the projection only once per time step, acknowledge that some accuracy is lost each time the mesh is coarsened and rely on very small time steps to keep it acceptable~\cite[Sect.~5.2]{CGRV2019}.

Several local alternatives have been proposed to avoid the global system. One family works element by element. It includes the B\'ezier projection~\cite{TSETE15}, later extended to hierarchical B-splines and reformulated in terms of univariate operators in~\cite{D'AR2020} and combined with the coarsening of THB-splines in~\cite{DGT2026}, and the element-wise least-squares transfer operators of~\cite{HALK2018}. Being local to each element, these projectors produce multi-valued control coefficients that must be averaged to recover a single spline, which introduces some smoothing and perturbs the solution also outside, near the boundary of the un-refined region. Closer to our approach, Jiang and Dolbow~\cite{JD2015} restrict the $L^2$-projection to the coarsening region, fitting only the coarse functions whose support is contained in it, so that the field is left unchanged outside; the fit still requires quadrature, and no bound on the committed error is provided. A second family replaces the projection by the quasi-interpolants available in hierarchical spaces~\cite{SM2016,Speleers2017}; this is the route followed, for instance, in the goal-adaptive shell simulations of~\cite{VMMB2024}. Quasi-interpolation is very cheap, since every coarse coefficient is computed independently of the others, but for the same reason it ignores the interaction of the overlapping coarse functions and offers no control over the part of the domain that should remain untouched. Closest to the present work is the coefficient-based data reduction operator of~\cite{BGS2017}. There, each THB-spline coefficient is obtained from a small least-squares problem posed directly on the coefficients of the target spline, with no pointwise evaluation, and the discrepancy between the coefficients of the target spline and those of its reduced representation, measured in the maximum norm and scaled by the range of the data, drives the construction of the hierarchical mesh. The setting is nevertheless the opposite of ours: a hierarchical representation of a given tensor-product spline is built by successive refinement from the coarsest level, rather than an existing hierarchical representation being un-refined. The coarse coefficients are moreover computed independently of one another, one local problem per function, so that the operator reproduces the coarsest tensor-product space but is not shown to reproduce the hierarchical space itself, and no part of the target spline is preserved exactly. With the exception of the coefficient-based criterion just mentioned, none of these approaches provides a computable bound for the accuracy lost in the transfer: it is either estimated a priori, with constants that are not available in practice~\cite{DGT2026}, or monitored through heuristic indicators with empirically tuned constants~\cite{JD2015,HALK2018}. We also mention~\cite{LSTHG}, where hierarchical refinement and coarsening are applied to moving-interface phase-field problems, and~\cite{HBJ2022}, which addresses the un-refinement of THB-spline parameterizations---a geometric problem whose requirements, such as preventing self-intersections, are different from ours.
 
The approach followed in this paper is purely algebraic: instead of
a projection that requires quadrature or point evaluations, we invert, in the least-squares sense, the refinement matrix that relates the coarse and the fine representations. Reversing a subdivision rule by local least squares goes back to~\cite{BS2000}, and was further analysed in~\cite{BGS2006}; in the spline setting, sharp local
error indicators for univariate knot removal were obtained in~\cite{FGM24}, and
local left inverses of tensor-product subdivision matrices were constructed
in~\cite{AFG25}.
The hierarchical case is not a level-by-level application of the
tensor-product operator. Un-refining a set of cells reactivates coarse
functions whose supports extend beyond the marked region, keeps some of their
children and removes the others, and touches active functions of every finer
level; in the truncated basis, it also changes the surviving coarse functions themselves. Our contribution is to identify these interactions and to show that all
of them can be resolved by local computations.
 
The construction rests on a single observation: after a suitable reordering, the
refinement matrix that expresses the coarsened basis in terms of the current one splits as
$K=\mathrm{diag}(K_{\loc},I)$, where the identity block acts on the functions
that the coarsening leaves untouched and $K_{\loc}$ is a small block-triangular
matrix. The coefficients of the reactivated B-splines are obtained from one local least-squares problem, and those of the retained B-splines involved in their expansion by explicit corrections. The resulting coarsening operator is a left
inverse of $K$ (Corollary~\ref{cor:exactness}): the transfer is exact whenever the fine spline already belongs
to the coarse space, and in all cases the spline is left unchanged outside the
supports of the functions that the step removes, a neighbourhood of the
coarsened region. The
maximum norm of the residual of the local
system is a computable a posteriori bound for the $L^\infty$-error with no unknown constants
(Proposition~\ref{prop:linf_bound}), which makes tolerance-driven coarsening possible. On top of
this operator we build an adaptive multi-level algorithm with an algebraic
rollback mechanism that keeps the accumulated error below a prescribed tolerance
(Theorem~\ref{thm:adaptive_success}). The whole construction carries over to THB-splines, where truncation makes the corrections unnecessary and the same constant-free bound holds (Proposition~\ref{prop:thb_linf_bound}). The procedure requires neither numerical
quadrature nor global assembly. The price of locality is a moderate loss of accuracy with respect to a global $L^2$-projection, which may modify every coefficient of the spline.
 
The paper is organized as follows. Section~\ref{sec:hierarchical_bsplines}
collects the notation on hierarchical splines and their subdivision and
multi-level transfer matrices.
Section~\ref{sec:coarsening} contains the core of the method: the block
structure of the refinement matrix, the local coarsening problem and its
solution, the constant-free error bound, and the resulting algorithms. Section~\ref{sec:adaptive_coarsening}
develops an adaptive multi-level algorithm with guaranteed error control. Section~\ref{sec:thb} extends the method to THB-splines, for which the corrections are no longer needed and only the right-hand side of the local problem changes.
Section~\ref{sec:numerical_experiments} presents numerical experiments along an adaptive refinement path and behind a moving front.
 
\section{Preliminaries on hierarchical B-splines}
\label{sec:hierarchical_bsplines}

In this section we introduce some notation and recall the fundamental concepts of
hierarchical B-splines~\cite{Kraft, VGJS2011}. Consider a sequence of nested
$d$-variate tensor-product spline spaces~(cf.~\cite{DeBoor,Schumi}) defined on
the parametric domain $\Omega = [0,1]^d$ and determined by a polynomial degree
$p$ in each direction and by nested knot vectors:
\begin{equation*}
    \mathcal{S}_0 \subset \mathcal{S}_1 \subset \dots \subset \mathcal{S}_\ell\subset \dots,
\end{equation*}
where each space $\mathcal{S}_\ell$ is spanned by the corresponding B-spline basis $\mathcal{B}_\ell := \{\beta_{i,\ell} \mid i = 1, \dots, N_\ell\}$, with $N_\ell:=\dim\mathcal S_\ell$. Let $\GG_\ell$ be the Cartesian mesh associated with $\mathcal{S}_\ell$. The cells $Q \in \GG_\ell$ are considered closed hyperrectangles with non-empty interiors. The support of a B-spline of level $\ell$ spans at most $p+1$ cells of $\GG_\ell$ in each direction, and exactly $(p+1)^d$ functions of $\mathcal{B}_\ell$ do not vanish identically on each cell of that level.
Throughout the paper the refinement between two consecutive levels is dyadic, so
that every cell of $\GG_\ell$ is the union of $2^d$ cells of $\GG_{\ell+1}$.

The algebraic foundation of hierarchical spline spaces is the two-scale relation: every B-spline of level $\ell$ is a linear combination of B-splines of level $\ell+1$,
\begin{equation}\label{E:two scale relation}
    \beta_{i,\ell} = \sum_{k=1}^{N_{\ell+1}} c_{k,\ell+1}(\beta_{i,\ell})\, \beta_{k,\ell+1}, \quad \forall \beta_{i,\ell} \in \mathcal{B}_\ell,
\end{equation}
with non-negative coefficients $c_{k,\ell+1}(\beta_{i,\ell}) \ge 0$. These coefficients form the \emph{subdivision matrix} $C_\ell^{\ell+1} \in \RR^{N_{\ell+1} \times N_\ell}$, with entries $(C_\ell^{\ell+1})_{ki} = c_{k,\ell+1}(\beta_{i,\ell})$; by the local support of B-splines, each column has only a few nonzero entries. We say that $\beta_{k,\ell+1}$ is a \emph{child} of $\beta_{i,\ell}$ if $c_{k,\ell+1}(\beta_{i,\ell}) > 0$, and denote by $\children{\beta_{i,\ell}}\subset\mathcal B_{\ell+1}$ the set of children of $\beta_{i,\ell}$ and by $\children{\mathcal F}$ the set of all the children of the functions in a set $\mathcal F$. The support of a child is contained in that of its parent, and a B-spline has at most $(p+2)^d$ children.

Similarly, a cell $Q'\in\GG_{\ell+1}$ is a child of $Q\in\GG_\ell$, and $Q$ is its parent, if $Q'\subset Q$. We denote by $\children{Q}$ the $2^d$ children of $Q$, and, for a set $\MM$ of cells, by $\children{\MM}$ and $\parents{\MM}$ the sets of all their children and of all their parents, respectively.

\paragraph{Hierarchical spline spaces.} A hierarchy of subdomains of $\Omega$ of depth $L \in \mathbb{N}$ is a sequence $\boldsymbol{\Omega}_L := \{\Omega_0, \dots, \Omega_L\}$ such that
\begin{equation*}
    \Omega = \Omega_0 \supseteq \Omega_1 \supseteq \cdots \supseteq \Omega_L = \varnothing,
\end{equation*}
where each $\Omega_{\ell}$, for $\ell \ge 1$, is a union of cells from $\GG_{\ell-1}$.

The \emph{hierarchical mesh} $\QQ$ is then defined by
\begin{equation*}
    \QQ := \bigcup_{\ell=0}^{L-1}\{ Q \in \GG_\ell \mid Q \subset \Omega_\ell \text{ and } Q \not\subset \Omega_{\ell+1} \},
\end{equation*}
and its elements are called \emph{active cells}. Thus $\GG_\ell$ always denotes the full Cartesian mesh of level $\ell$, whereas $\QQ\cap\GG_\ell$ collects its active cells. This mirrors the convention used below for basis functions, where $\mathcal{B}_\ell$ is the full tensor-product basis of level $\ell$ and $\mathcal{F}^A_\ell$ its active part.

The \emph{hierarchical basis} $\mathcal{H}$ is defined recursively starting from $\mathcal{H}_0 := \mathcal{B}_0$. For $\ell = 0, \ldots, L-2$, $\mathcal H_{\ell+1}$ is obtained from $\mathcal H_\ell$ by replacing the functions of $\mathcal H_\ell$ whose support is contained in $\Omega_{\ell+1}$ by all the B-splines of level $\ell+1$ with the same property:
\begin{equation*}
    \mathcal{H}_{\ell+1} := \{ \beta \in \mathcal{H}_\ell \mid \operatorname{supp}(\beta) \not\subset \Omega_{\ell+1} \} \cup \{ \beta \in \mathcal{B}_{\ell+1} \mid \operatorname{supp}(\beta) \subset \Omega_{\ell+1} \}.
\end{equation*}
The hierarchical basis is then defined by $\mathcal{H} := \mathcal{H}_{L-1}$, and the corresponding hierarchical spline space is $\operatorname{span}\mathcal{H}$.

This recursive construction partitions each single-level basis $\mathcal{B}_\ell$ into three mutually exclusive sets:
\begin{itemize}
    \item \textbf{Active functions} ($\mathcal{F}^A_\ell := \mathcal{H} \cap \mathcal{B}_\ell$): B-splines retained in the hierarchical basis $\mathcal{H}$. Their support contains at least one active cell of level $\ell$.
    \item  \textbf{Deactivated functions} ($\mathcal{F}^D_\ell := \{ \beta \in \mathcal{B}_\ell \mid \operatorname{supp}(\beta) \subset \Omega_{\ell+1} \}$): B-splines that enter the construction at level $\ell$ but are discarded at level $\ell+1$, because their support is contained in the refined subdomain $\Omega_{\ell+1}$.
    \item \textbf{Passive functions}: B-splines whose support is not contained in $\Omega_\ell$, which never participate in the construction.
\end{itemize}

We denote the cardinalities of these sets by $N_\ell^A := \#\mathcal{F}^A_\ell$ and $N_\ell^D := \#\mathcal{F}^D_\ell$. The global set of active functions is $\mathcal{F}^A := \bigcup_{\ell=0}^{L-1} \mathcal{F}^A_\ell = \mathcal{H}$, and the global set of deactivated functions is $\mathcal{F}^D := \bigcup_{\ell=0}^{L-1} \mathcal{F}^D_\ell$.
A spline $s\in\Span\HH$ is written as
$s=\sum_{\beta\in\mathcal F^A}c_\beta\,\beta$, and we call $c\in\RR^{n}$,
with $n:=\#\mathcal F^A$, its coefficient vector.

\paragraph{Notation.} We address submatrices and subvectors by sets of B-splines rather than by integers: for a matrix $M$ and a vector $v$, we write $M(\mathcal{X}, \mathcal{Y})$ for the submatrix with rows indexed by $\mathcal{X}$ and columns indexed by $\mathcal{Y}$, and $v(\mathcal{X})$ for the subvector with components indexed by $\mathcal{X}$. Two canonical identifications are involved, and the object being addressed determines which one applies: every B-spline of level $\ell$ has a unique position in the tensor-product basis $\mathcal{B}_\ell$, which is the one used for level-wise matrices such as $C_\ell^{\ell+1}$; an active B-spline has a unique position in the hierarchical basis $\mathcal{H}$, which is the one used for coefficient vectors such as $c$. A matrix may use different identifications for its rows and its columns, as the transfer matrices $M_\ell$ below. No ambiguity arises, and we therefore keep the local-to-global maps implicit.

\paragraph{Multi-level transfer matrices for deactivated B-splines.} The coarsening step requires each deactivated B-spline of level $\ell$ to be expressed in terms of the active functions of the finer levels. Between two consecutive levels, this information is contained in the blocks
\begin{equation*}
    C_{\ell}^{\ell+1, A} := C_{\ell}^{\ell+1}(\mathcal{F}_{\ell+1}^A, \mathcal{F}_\ell^D)\in\RR^{N_{\ell+1}^A \times N_\ell^D}
    \qquad \text{and} \qquad
    C_{\ell}^{\ell+1, D} := C_{\ell}^{\ell+1}(\mathcal{F}_{\ell+1}^D, \mathcal{F}_\ell^D)\in\RR^{N_{\ell+1}^D \times N_\ell^D},
\end{equation*}
which contain the coefficients of the deactivated functions of level $\ell$ along the active and the deactivated functions of level $\ell+1$, respectively. No other rows are needed, since the children of a deactivated function have their supports in $\Omega_{\ell+1}$ and are therefore active or deactivated.

Composing these blocks across the levels, we define for $\ell=0,\dots,L-2$ the \emph{multi-level transfer matrix} $M_\ell\in\RR^{N^A_{>\ell}\times N^D_\ell}$, where $\mathcal F^A_{>\ell}:=\bigcup_{k=\ell+1}^{L-1}\mathcal F^A_k$ collects the active functions of the finer levels and $N^A_{>\ell}:=\#\mathcal F^A_{>\ell}$, by the recursion
\begin{equation}\label{E:multilevel transfer matrices}
    M_{L-2}:=C_{L-2}^{L-1,A},
    \qquad
    M_\ell:=\begin{pmatrix} C_{\ell}^{\ell+1,A}\\[4pt] M_{\ell+1}\,C_{\ell}^{\ell+1,D}\end{pmatrix}
    \quad (0\le\ell\le L-3),
\end{equation}
the initialization being consistent because $\mathcal F^D_{L-1}=\varnothing$. The first block accounts for the active children of level $\ell+1$, and the second one propagates the contribution of the deactivated children to their own active descendants. Consequently, the column of $M_\ell$ indexed by $\beta\in\mathcal F^D_\ell$ contains the coefficients of $\beta$ in terms of the active functions of the levels $\ell+1,\dots,L-1$. Its columns are indexed by $\mathcal F^D_\ell$ within $\mathcal B_\ell$, and its rows by $\mathcal F^A_{>\ell}$ in the order of the hierarchical basis. These matrices depend only on the hierarchical space, and can be assembled from the univariate subdivision matrices stored by GeoPDEs~\cite{geopdes,GV18}.
 
\section{Coarsening for hierarchical B-splines}
\label{sec:coarsening}
 
Refining a hierarchical mesh changes the basis, not the spline: the coefficients in the fine basis are obtained from the coarse ones by a refinement matrix. Coarsening is the inverse operation, and it is exact only for splines that already belong to the coarse space. The procedure proposed in this section realizes a left inverse of the refinement matrix, without ever assembling it, and returns for every other spline a computable bound on what is lost. Section~\ref{S:local_block_structure} identifies the basis functions involved and shows that the refinement matrix splits into a small local block $K_{\loc}$ and an identity block; Section~\ref{S:local_coarsening_problem} formulates the local coarsening problem, solves it in closed form, shows that the resulting operator is a left inverse of the refinement matrix and a projection onto the coarse space, and derives a constant-free error bound; Section~\ref{sec:algorithmic_implementation} gives the algorithms.

 \subsection{From refinement to coarsening: the local block structure}
\label{S:local_block_structure}

Let $(\hat\QQ,\hat\HH)$ be a hierarchical mesh and its hierarchical basis, with hierarchy of subdomains $\hat{\boldsymbol\Omega}_L=\{\hat\Omega_k\}_{k=0}^{L}$, and let $\hat\QQ_{\loc}\subset\hat\QQ\cap\GG_\ell$ be a set of active cells of level $\ell$, with $\ell\le L-2$ (if $\ell=L-1$, the depth $L$ is first increased by one). Refining these cells produces the hierarchical mesh $\QQ$ and basis $\HH$ associated with the hierarchy of subdomains $\boldsymbol\Omega_L=\{\Omega_k\}_{k=0}^{L}$ given by 
\begin{equation}\label{eq:refined_hierarchy}
    \Omega_{\ell+1}:=\hat\Omega_{\ell+1}\cup\Omega_{\loc},
    \qquad
    \Omega_k:=\hat\Omega_k\quad\text{for }k\neq\ell+1,
\end{equation}
where $\QQ_{\loc}:=\children{\hat\QQ_{\loc}}\subset\GG_{\ell+1}$ are the new fine cells and
\begin{equation}\label{eq:local_region}
    \Omega_{\loc}:=\bigcup_{Q\in\hat\QQ_{\loc}}Q=\bigcup_{Q'\in\QQ_{\loc}}Q'
\end{equation}
is the region where the mesh changes. Throughout, $|\cdot|$ denotes the Lebesgue measure in $\mathbb{R}^d$, and we use conditions of the form $|\supp\beta\cap\Omega_{\loc}|>0$ rather than $\supp\beta\cap\Omega_{\loc}\neq\varnothing$, since both supports and cells are closed sets and functions merely touching $\partial\Omega_{\loc}$ must not be regarded as interacting with the local region. For brevity, we say that two subsets of $\Omega$ \emph{overlap} if their intersection has positive measure, and that a function overlaps a set if its support does. Since $\Span\hat\HH\subseteq\Span\HH$, each spline $s\in\Span\hat\HH$ can also be represented in the fine basis through a refinement matrix $K\in\RR^{n\times\hat n}$, where $\hat n:=\#\hat\HH$ and $n:=\#\HH$. More precisely, if
\[
    s=\sum_{\hat\beta\in\hat\HH}\hat c_{\hat\beta}\hat\beta
    =\sum_{\beta\in\HH}c_\beta\beta,
\]
then the two coefficient vectors are related by
\begin{equation}\label{eq:refinement_relation}
    c=K\hat c.
\end{equation}
Thus, the matrix $K$ contains the purely algebraic information needed to pass from the coarse representation to the fine one. The nonzero entries of a column of $K$ are precisely the coefficients of the corresponding coarse basis function in the fine basis. In particular, the column of a coarse function that remains active after the refinement is a unit vector, since that function also belongs to the fine basis, whereas the column of a function deactivated by the refinement contains the coefficients of its expansion in the fine basis.
 
Only a few basis functions, all of them near $\Omega_{\loc}$, change their status or their coefficient; they are identified precisely below. This is what makes the refinement matrix useful for coarsening: it separates the coefficient vector into a small part that changes and a large part that does not.

\paragraph{Classification of the basis functions.} We now identify the B-splines whose status differs between the two configurations, and those that are active in both but whose coefficient is in general modified when a spline is rewritten from the coarse basis into the fine one. All the sets below are defined \emph{in terms of the fine configuration} $(\QQ,\HH)$ and its hierarchy $\boldsymbol\Omega_L$, which is the one available when coarsening is performed; their status in both configurations and the treatment of their coefficients by the coarsening step are summarized in Table~\ref{tab:local_sets}. 
 
\begin{table}[htb]
    \small
    \centering
    \renewcommand{\arraystretch}{1.25}
    \begin{tabular}{llccc}
        \toprule
        \textbf{Set} & \textbf{Level} & \textbf{Status in $\hat\HH$} & \textbf{Status in $\HH$} & \textbf{Coefficient} \\
        \midrule
        $\mathcal A_{\loc}$ & $\ell$ & active & deactivated & recovered by least squares \\
        $\mathcal C_{\loc}^{\fineonly}$ & $\ell+1$ & passive & active & removed \\
        $\mathcal C_{\loc}^{\mathrm{coarse}}$ & $\ell+1$ & active & active & corrected \\
        $\mathcal C_{\loc}^{\mathrm{deact}}$ & $\ell+1$ & deactivated & deactivated & --- \\
        $\mathcal D_{\loc}$ & $\ge\ell+2$ & active & active & corrected \\
        $\mathcal H_{\mathrm{rem}}$ & any & active & active & copied \\
        \bottomrule
    \end{tabular}
    \caption{The local B-spline sets, their status in the coarse basis $\hat\HH$ and in the fine basis $\HH$, and the way their coefficients are treated by the coarsening procedure.}
    \label{tab:local_sets}
\end{table}

The \emph{ancestors} are the level-$\ell$ B-splines that are deactivated in the fine configuration and are reactivated by the coarsening,
\begin{equation}\label{eq:def_Aloc}
    \mathcal A_{\loc}:=\{\beta\in\mathcal F_\ell^D:\ |\supp\beta\cap\Omega_{\loc}|>0\}.
\end{equation}
Their children $\mathcal C_{\loc}:=\children{\mathcal A_{\loc}}$ satisfy $\supp\gamma\subset\supp\beta\subset\Omega_{\ell+1}$ for every $\gamma\in\mathcal C_{\loc}$, where $\beta\in\mathcal A_{\loc}$ is a parent of $\gamma$, so that each of them is either active or deactivated at level $\ell+1$. This yields the disjoint splitting
\begin{equation}\label{eq:children_split}
    \mathcal C_{\loc}^{\mathrm{fine}}:=\mathcal C_{\loc}\cap\mathcal F_{\ell+1}^A,
    \qquad
    \mathcal C_{\loc}^{\mathrm{deact}}:=\mathcal C_{\loc}\cap\mathcal F_{\ell+1}^D.
\end{equation}
The fine basis functions that are removed by the coarsening step are
\begin{equation}\label{eq:def_fineonly}
    \mathcal C_{\loc}^{\fineonly}:=\{\beta\in\mathcal F_{\ell+1}^A:\ |\supp\beta\cap\Omega_{\loc}|>0\},
\end{equation}
whereas the active children of the ancestors that are not removed, and are therefore active in both configurations, are
\begin{equation}\label{eq:def_coarse_children}
    \mathcal C_{\loc}^{\mathrm{coarse}}:=\mathcal C_{\loc}^{\mathrm{fine}}\setminus\mathcal C_{\loc}^{\fineonly}
    =\{\beta\in\mathcal C_{\loc}^{\mathrm{fine}}:\ |\supp\beta\cap\Omega_{\loc}|=0\}.
\end{equation}
 
Definition~\eqref{eq:def_fineonly} scans the active functions of level $\ell+1$, rather than the children of the ancestors only, as in the alternative $\{\beta\in\mathcal C_{\loc}^{\mathrm{fine}}:|\supp\beta\cap\Omega_{\loc}|>0\}$. The latter set is contained in the former, and the inclusion can be strict, since the children of the deactivated B-splines of level $\ell$ need not exhaust the B-splines activated at level $\ell+1$; see~\cite[Sect.~5 and Fig.~4]{BG17}. Searching among the children would then leave out active B-splines that the coarsening removes nonetheless; see Remark~\ref{rem:no_ancestors} below.

Finally, if $\mathcal C_{\loc}^{\mathrm{deact}}\neq\varnothing$, that is, some children of the ancestors are themselves deactivated, the expansion of the ancestors in the fine basis involves the active descendants of these children, which are reached through the multi-level transfer matrix $M_{\ell+1}$ given by~\eqref{E:multilevel transfer matrices}. The set $\mathcal D_{\loc}$ of active descendants involved and the associated local multi-level transfer matrix $M_{\loc}$ are given by
\begin{equation}\label{eq:def_Dloc}
    \mathcal D_{\loc}:=\Big\{\gamma\in\textstyle\bigcup_{k=\ell+2}^{L-1}\mathcal F_k^A:\
    M_{\ell+1}(\{\gamma\},\mathcal C_{\loc}^{\mathrm{deact}})\neq\mathbf 0\Big\},
    \qquad
    M_{\loc}:=M_{\ell+1}(\mathcal D_{\loc},\mathcal C_{\loc}^{\mathrm{deact}}),
\end{equation}
and $\mathcal H_{\mathrm{rem}}:=\HH\setminus(\mathcal C_{\loc}^{\fineonly}\cup\mathcal C_{\loc}^{\mathrm{coarse}}\cup\mathcal D_{\loc})$ denotes all the remaining active functions, which belong to both bases and have the same coefficient in both representations. For the finest transition $\ell=L-2$ one has $\mathcal F_{L-1}^D=\varnothing$, hence $\mathcal C_{\loc}^{\mathrm{deact}}=\varnothing$ and $\mathcal D_{\loc}=\varnothing$, and the matrices $M_{\loc}$ and $K_D$ below are empty. The matrix $M_{\loc}$ contains the coefficients that express the deactivated children in terms of their active descendants; hence, even when the hierarchy has more than two levels, the expansion of an ancestor in the fine basis is obtained from its two-scale relation~\eqref{E:two scale relation}, followed by $M_{\loc}$ for its deactivated children. 
 
Applying the two-scale relation~\eqref{E:two scale relation} successively, one obtains the action of $K$, that is, how the coefficients of a spline change when it is rewritten from the coarse basis into the fine one: the coefficients of $\mathcal A_{\loc}$ determine those of $\mathcal C_{\loc}^{\fineonly}$ and add a term to those of $\mathcal C_{\loc}^{\mathrm{coarse}}$ and $\mathcal D_{\loc}$, whereas the coefficients of $\mathcal H_{\mathrm{rem}}$ are copied unchanged. Figure~\ref{fig:coarsening_sets} displays these sets on a concrete hierarchical mesh, and in the block form of $K$ below, the three actions just described are precisely the three nonzero blocks of its first block column.
 
\begin{figure}[htbp]
    \centering
    \definecolor{setA}{RGB}{0,114,178}
    \definecolor{setF}{RGB}{213,94,0}
    \definecolor{setC}{RGB}{0,158,115}
    \definecolor{setD}{RGB}{204,121,167}
    \definecolor{markedbg}{RGB}{245,222,205}
    \begin{tikzpicture}[
        x=1cm, y=1cm,
        cell/.style={draw=black!40, line width=0.3pt},
        dombox/.style={draw=black!75, line width=0.7pt},
        markedfill/.style={fill=markedbg},
        markedline/.style={draw=setF, line width=1.3pt},
        ghostline/.style={draw=setF!55, line width=1.0pt, dotted},
        suppline/.style={draw=setA, line width=0.9pt, dashed, dash pattern=on 3pt off 2pt},
        mkAon/.style={circle, draw=setA, fill=setA, inner sep=0pt, minimum size=6.4pt},
        mkAoff/.style={circle, draw=setA, fill=white, line width=1.0pt,
                       inner sep=0pt, minimum size=6.4pt},
        mkF/.style={rectangle, draw=setF, fill=setF, inner sep=0pt, minimum size=5.0pt},
        mkC/.style={diamond, draw=setC, fill=setC, inner sep=0pt, minimum size=5.4pt},
        mkD/.style={circle, draw=setD, fill=setD, inner sep=0pt, minimum size=3.0pt},
        panellabel/.style={font=\small, anchor=north},
        leg/.style={font=\footnotesize, anchor=west},
    ]
    \begin{scope}[shift={(0.0,0)}]
      \fill[markedfill] (1,1) rectangle (3,3);
      \draw[cell] (0,0) rectangle (1,1);
      \draw[cell] (1,0) rectangle (2,1);
      \draw[cell] (2,0) rectangle (3,1);
      \draw[cell] (3,0) rectangle (4,1);
      \draw[cell] (4,0) rectangle (5,1);
      \draw[cell] (5,0) rectangle (6,1);
      \draw[cell] (0,1) rectangle (1,2);
      \draw[cell] (5,1) rectangle (6,2);
      \draw[cell] (0,2) rectangle (1,3);
      \draw[cell] (5,2) rectangle (6,3);
      \draw[cell] (0,3) rectangle (1,4);
      \draw[cell] (5,3) rectangle (6,4);
      \draw[cell] (0,4) rectangle (1,5);
      \draw[cell] (5,4) rectangle (6,5);
      \draw[cell] (0,5) rectangle (1,6);
      \draw[cell] (1,5) rectangle (2,6);
      \draw[cell] (2,5) rectangle (3,6);
      \draw[cell] (3,5) rectangle (4,6);
      \draw[cell] (4,5) rectangle (5,6);
      \draw[cell] (5,5) rectangle (6,6);
      \draw[cell] (1,1) rectangle (1.5,1.5);
      \draw[cell] (1.5,1) rectangle (2,1.5);
      \draw[cell] (2,1) rectangle (2.5,1.5);
      \draw[cell] (2.5,1) rectangle (3,1.5);
      \draw[cell] (3,1) rectangle (3.5,1.5);
      \draw[cell] (3.5,1) rectangle (4,1.5);
      \draw[cell] (4,1) rectangle (4.5,1.5);
      \draw[cell] (4.5,1) rectangle (5,1.5);
      \draw[cell] (1,1.5) rectangle (1.5,2);
      \draw[cell] (1.5,1.5) rectangle (2,2);
      \draw[cell] (2,1.5) rectangle (2.5,2);
      \draw[cell] (2.5,1.5) rectangle (3,2);
      \draw[cell] (3,1.5) rectangle (3.5,2);
      \draw[cell] (3.5,1.5) rectangle (4,2);
      \draw[cell] (4,1.5) rectangle (4.5,2);
      \draw[cell] (4.5,1.5) rectangle (5,2);
      \draw[cell] (1,2) rectangle (1.5,2.5);
      \draw[cell] (1.5,2) rectangle (2,2.5);
      \draw[cell] (2,2) rectangle (2.5,2.5);
      \draw[cell] (2.5,2) rectangle (3,2.5);
      \draw[cell] (3,2) rectangle (3.5,2.5);
      \draw[cell] (3.5,2) rectangle (4,2.5);
      \draw[cell] (4,2) rectangle (4.5,2.5);
      \draw[cell] (4.5,2) rectangle (5,2.5);
      \draw[cell] (1,2.5) rectangle (1.5,3);
      \draw[cell] (1.5,2.5) rectangle (2,3);
      \draw[cell] (2,2.5) rectangle (2.5,3);
      \draw[cell] (2.5,2.5) rectangle (3,3);
      \draw[cell] (3,2.5) rectangle (3.5,3);
      \draw[cell] (3.5,2.5) rectangle (4,3);
      \draw[cell] (4,2.5) rectangle (4.5,3);
      \draw[cell] (4.5,2.5) rectangle (5,3);
      \draw[cell] (1,3) rectangle (1.5,3.5);
      \draw[cell] (1.5,3) rectangle (2,3.5);
      \draw[cell] (2,3) rectangle (2.5,3.5);
      \draw[cell] (2.5,3) rectangle (3,3.5);
      \draw[cell] (1,3.5) rectangle (1.5,4);
      \draw[cell] (1.5,3.5) rectangle (2,4);
      \draw[cell] (2,3.5) rectangle (2.5,4);
      \draw[cell] (2.5,3.5) rectangle (3,4);
      \draw[cell] (1,4) rectangle (1.5,4.5);
      \draw[cell] (1.5,4) rectangle (2,4.5);
      \draw[cell] (2,4) rectangle (2.5,4.5);
      \draw[cell] (2.5,4) rectangle (3,4.5);
      \draw[cell] (1,4.5) rectangle (1.5,5);
      \draw[cell] (1.5,4.5) rectangle (2,5);
      \draw[cell] (2,4.5) rectangle (2.5,5);
      \draw[cell] (2.5,4.5) rectangle (3,5);
      \draw[cell] (3,3) rectangle (3.25,3.25);
      \draw[cell] (3.25,3) rectangle (3.5,3.25);
      \draw[cell] (3.5,3) rectangle (3.75,3.25);
      \draw[cell] (3.75,3) rectangle (4,3.25);
      \draw[cell] (4,3) rectangle (4.25,3.25);
      \draw[cell] (4.25,3) rectangle (4.5,3.25);
      \draw[cell] (4.5,3) rectangle (4.75,3.25);
      \draw[cell] (4.75,3) rectangle (5,3.25);
      \draw[cell] (3,3.25) rectangle (3.25,3.5);
      \draw[cell] (3.25,3.25) rectangle (3.5,3.5);
      \draw[cell] (3.5,3.25) rectangle (3.75,3.5);
      \draw[cell] (3.75,3.25) rectangle (4,3.5);
      \draw[cell] (4,3.25) rectangle (4.25,3.5);
      \draw[cell] (4.25,3.25) rectangle (4.5,3.5);
      \draw[cell] (4.5,3.25) rectangle (4.75,3.5);
      \draw[cell] (4.75,3.25) rectangle (5,3.5);
      \draw[cell] (3,3.5) rectangle (3.25,3.75);
      \draw[cell] (3.25,3.5) rectangle (3.5,3.75);
      \draw[cell] (3.5,3.5) rectangle (3.75,3.75);
      \draw[cell] (3.75,3.5) rectangle (4,3.75);
      \draw[cell] (4,3.5) rectangle (4.25,3.75);
      \draw[cell] (4.25,3.5) rectangle (4.5,3.75);
      \draw[cell] (4.5,3.5) rectangle (4.75,3.75);
      \draw[cell] (4.75,3.5) rectangle (5,3.75);
      \draw[cell] (3,3.75) rectangle (3.25,4);
      \draw[cell] (3.25,3.75) rectangle (3.5,4);
      \draw[cell] (3.5,3.75) rectangle (3.75,4);
      \draw[cell] (3.75,3.75) rectangle (4,4);
      \draw[cell] (4,3.75) rectangle (4.25,4);
      \draw[cell] (4.25,3.75) rectangle (4.5,4);
      \draw[cell] (4.5,3.75) rectangle (4.75,4);
      \draw[cell] (4.75,3.75) rectangle (5,4);
      \draw[cell] (3,4) rectangle (3.25,4.25);
      \draw[cell] (3.25,4) rectangle (3.5,4.25);
      \draw[cell] (3.5,4) rectangle (3.75,4.25);
      \draw[cell] (3.75,4) rectangle (4,4.25);
      \draw[cell] (4,4) rectangle (4.25,4.25);
      \draw[cell] (4.25,4) rectangle (4.5,4.25);
      \draw[cell] (4.5,4) rectangle (4.75,4.25);
      \draw[cell] (4.75,4) rectangle (5,4.25);
      \draw[cell] (3,4.25) rectangle (3.25,4.5);
      \draw[cell] (3.25,4.25) rectangle (3.5,4.5);
      \draw[cell] (3.5,4.25) rectangle (3.75,4.5);
      \draw[cell] (3.75,4.25) rectangle (4,4.5);
      \draw[cell] (4,4.25) rectangle (4.25,4.5);
      \draw[cell] (4.25,4.25) rectangle (4.5,4.5);
      \draw[cell] (4.5,4.25) rectangle (4.75,4.5);
      \draw[cell] (4.75,4.25) rectangle (5,4.5);
      \draw[cell] (3,4.5) rectangle (3.25,4.75);
      \draw[cell] (3.25,4.5) rectangle (3.5,4.75);
      \draw[cell] (3.5,4.5) rectangle (3.75,4.75);
      \draw[cell] (3.75,4.5) rectangle (4,4.75);
      \draw[cell] (4,4.5) rectangle (4.25,4.75);
      \draw[cell] (4.25,4.5) rectangle (4.5,4.75);
      \draw[cell] (4.5,4.5) rectangle (4.75,4.75);
      \draw[cell] (4.75,4.5) rectangle (5,4.75);
      \draw[cell] (3,4.75) rectangle (3.25,5);
      \draw[cell] (3.25,4.75) rectangle (3.5,5);
      \draw[cell] (3.5,4.75) rectangle (3.75,5);
      \draw[cell] (3.75,4.75) rectangle (4,5);
      \draw[cell] (4,4.75) rectangle (4.25,5);
      \draw[cell] (4.25,4.75) rectangle (4.5,5);
      \draw[cell] (4.5,4.75) rectangle (4.75,5);
      \draw[cell] (4.75,4.75) rectangle (5,5);
      \draw[dombox] (0,0) rectangle (6,6);
      \draw[markedline] (1,1) rectangle (3,3);
      \draw[suppline] (2,2) rectangle (5,5); 
      \node[mkD] at (3.375,3.375) {};
      \node[mkD] at (3.625,3.375) {};
      \node[mkD] at (3.875,3.375) {};
      \node[mkD] at (4.125,3.375) {};
      \node[mkD] at (4.375,3.375) {};
      \node[mkD] at (4.625,3.375) {};
      \node[mkD] at (3.375,3.625) {};
      \node[mkD] at (3.625,3.625) {};
      \node[mkD] at (3.875,3.625) {};
      \node[mkD] at (4.125,3.625) {};
      \node[mkD] at (4.375,3.625) {};
      \node[mkD] at (4.625,3.625) {};
      \node[mkD] at (3.375,3.875) {};
      \node[mkD] at (3.625,3.875) {};
      \node[mkD] at (3.875,3.875) {};
      \node[mkD] at (4.125,3.875) {};
      \node[mkD] at (4.375,3.875) {};
      \node[mkD] at (4.625,3.875) {};
      \node[mkD] at (3.375,4.125) {};
      \node[mkD] at (3.625,4.125) {};
      \node[mkD] at (3.875,4.125) {};
      \node[mkD] at (4.125,4.125) {};
      \node[mkD] at (4.375,4.125) {};
      \node[mkD] at (4.625,4.125) {};
      \node[mkD] at (3.375,4.375) {};
      \node[mkD] at (3.625,4.375) {};
      \node[mkD] at (3.875,4.375) {};
      \node[mkD] at (4.125,4.375) {};
      \node[mkD] at (4.375,4.375) {};
      \node[mkD] at (4.625,4.375) {};
      \node[mkD] at (3.375,4.625) {};
      \node[mkD] at (3.625,4.625) {};
      \node[mkD] at (3.875,4.625) {};
      \node[mkD] at (4.125,4.625) {};
      \node[mkD] at (4.375,4.625) {};
      \node[mkD] at (4.625,4.625) {};
      \node[mkC] at (1.75,3.75) {};
      \node[mkC] at (1.75,4.25) {};
      \node[mkC] at (2.25,3.75) {};
      \node[mkC] at (2.25,4.25) {};
      \node[mkC] at (2.75,3.75) {};
      \node[mkC] at (2.75,4.25) {};
      \node[mkC] at (3.25,3.75) {};
      \node[mkC] at (3.25,4.25) {};
      \node[mkC] at (3.75,1.75) {};
      \node[mkC] at (3.75,2.25) {};
      \node[mkC] at (3.75,2.75) {};
      \node[mkC] at (3.75,3.25) {};
      \node[mkC] at (4.25,1.75) {};
      \node[mkC] at (4.25,2.25) {};
      \node[mkC] at (4.25,2.75) {};
      \node[mkC] at (4.25,3.25) {};
      \node[mkF] at (1.75,1.75) {};
      \node[mkF] at (2.25,1.75) {};
      \node[mkF] at (2.75,1.75) {};
      \node[mkF] at (3.25,1.75) {};
      \node[mkF] at (1.75,2.25) {};
      \node[mkF] at (2.25,2.25) {};
      \node[mkF] at (2.75,2.25) {};
      \node[mkF] at (3.25,2.25) {};
      \node[mkF] at (1.75,2.75) {};
      \node[mkF] at (2.25,2.75) {};
      \node[mkF] at (2.75,2.75) {};
      \node[mkF] at (3.25,2.75) {};
      \node[mkF] at (1.75,3.25) {};
      \node[mkF] at (2.25,3.25) {};
      \node[mkF] at (2.75,3.25) {};
      \node[mkF] at (3.25,3.25) {};
      \node[mkAoff] at (2.5,2.5) {};
      \node[mkAoff] at (3.5,2.5) {};
      \node[mkAoff] at (2.5,3.5) {};
      \node[mkAoff] at (3.5,3.5) {};
      \node[panellabel] at (3,-0.45) {(a) fine configuration $\mathcal{Q}$, $\mathcal{H}$};
    \end{scope}
 
    \begin{scope}[shift={(7.3,0)}]
      \draw[cell] (0,0) rectangle (1,1);
      \draw[cell] (1,0) rectangle (2,1);
      \draw[cell] (2,0) rectangle (3,1);
      \draw[cell] (3,0) rectangle (4,1);
      \draw[cell] (4,0) rectangle (5,1);
      \draw[cell] (5,0) rectangle (6,1);
      \draw[cell] (0,1) rectangle (1,2);
      \draw[cell] (1,1) rectangle (2,2);
      \draw[cell] (2,1) rectangle (3,2);
      \draw[cell] (5,1) rectangle (6,2);
      \draw[cell] (0,2) rectangle (1,3);
      \draw[cell] (1,2) rectangle (2,3);
      \draw[cell] (2,2) rectangle (3,3);
      \draw[cell] (5,2) rectangle (6,3);
      \draw[cell] (0,3) rectangle (1,4);
      \draw[cell] (5,3) rectangle (6,4);
      \draw[cell] (0,4) rectangle (1,5);
      \draw[cell] (5,4) rectangle (6,5);
      \draw[cell] (0,5) rectangle (1,6);
      \draw[cell] (1,5) rectangle (2,6);
      \draw[cell] (2,5) rectangle (3,6);
      \draw[cell] (3,5) rectangle (4,6);
      \draw[cell] (4,5) rectangle (5,6);
      \draw[cell] (5,5) rectangle (6,6);
      \draw[cell] (3,1) rectangle (3.5,1.5);
      \draw[cell] (3.5,1) rectangle (4,1.5);
      \draw[cell] (4,1) rectangle (4.5,1.5);
      \draw[cell] (4.5,1) rectangle (5,1.5);
      \draw[cell] (3,1.5) rectangle (3.5,2);
      \draw[cell] (3.5,1.5) rectangle (4,2);
      \draw[cell] (4,1.5) rectangle (4.5,2);
      \draw[cell] (4.5,1.5) rectangle (5,2);
      \draw[cell] (3,2) rectangle (3.5,2.5);
      \draw[cell] (3.5,2) rectangle (4,2.5);
      \draw[cell] (4,2) rectangle (4.5,2.5);
      \draw[cell] (4.5,2) rectangle (5,2.5);
      \draw[cell] (3,2.5) rectangle (3.5,3);
      \draw[cell] (3.5,2.5) rectangle (4,3);
      \draw[cell] (4,2.5) rectangle (4.5,3);
      \draw[cell] (4.5,2.5) rectangle (5,3);
      \draw[cell] (1,3) rectangle (1.5,3.5);
      \draw[cell] (1.5,3) rectangle (2,3.5);
      \draw[cell] (2,3) rectangle (2.5,3.5);
      \draw[cell] (2.5,3) rectangle (3,3.5);
      \draw[cell] (1,3.5) rectangle (1.5,4);
      \draw[cell] (1.5,3.5) rectangle (2,4);
      \draw[cell] (2,3.5) rectangle (2.5,4);
      \draw[cell] (2.5,3.5) rectangle (3,4);
      \draw[cell] (1,4) rectangle (1.5,4.5);
      \draw[cell] (1.5,4) rectangle (2,4.5);
      \draw[cell] (2,4) rectangle (2.5,4.5);
      \draw[cell] (2.5,4) rectangle (3,4.5);
      \draw[cell] (1,4.5) rectangle (1.5,5);
      \draw[cell] (1.5,4.5) rectangle (2,5);
      \draw[cell] (2,4.5) rectangle (2.5,5);
      \draw[cell] (2.5,4.5) rectangle (3,5);
      \draw[cell] (3,3) rectangle (3.25,3.25);
      \draw[cell] (3.25,3) rectangle (3.5,3.25);
      \draw[cell] (3.5,3) rectangle (3.75,3.25);
      \draw[cell] (3.75,3) rectangle (4,3.25);
      \draw[cell] (4,3) rectangle (4.25,3.25);
      \draw[cell] (4.25,3) rectangle (4.5,3.25);
      \draw[cell] (4.5,3) rectangle (4.75,3.25);
      \draw[cell] (4.75,3) rectangle (5,3.25);
      \draw[cell] (3,3.25) rectangle (3.25,3.5);
      \draw[cell] (3.25,3.25) rectangle (3.5,3.5);
      \draw[cell] (3.5,3.25) rectangle (3.75,3.5);
      \draw[cell] (3.75,3.25) rectangle (4,3.5);
      \draw[cell] (4,3.25) rectangle (4.25,3.5);
      \draw[cell] (4.25,3.25) rectangle (4.5,3.5);
      \draw[cell] (4.5,3.25) rectangle (4.75,3.5);
      \draw[cell] (4.75,3.25) rectangle (5,3.5);
      \draw[cell] (3,3.5) rectangle (3.25,3.75);
      \draw[cell] (3.25,3.5) rectangle (3.5,3.75);
      \draw[cell] (3.5,3.5) rectangle (3.75,3.75);
      \draw[cell] (3.75,3.5) rectangle (4,3.75);
      \draw[cell] (4,3.5) rectangle (4.25,3.75);
      \draw[cell] (4.25,3.5) rectangle (4.5,3.75);
      \draw[cell] (4.5,3.5) rectangle (4.75,3.75);
      \draw[cell] (4.75,3.5) rectangle (5,3.75);
      \draw[cell] (3,3.75) rectangle (3.25,4);
      \draw[cell] (3.25,3.75) rectangle (3.5,4);
      \draw[cell] (3.5,3.75) rectangle (3.75,4);
      \draw[cell] (3.75,3.75) rectangle (4,4);
      \draw[cell] (4,3.75) rectangle (4.25,4);
      \draw[cell] (4.25,3.75) rectangle (4.5,4);
      \draw[cell] (4.5,3.75) rectangle (4.75,4);
      \draw[cell] (4.75,3.75) rectangle (5,4);
      \draw[cell] (3,4) rectangle (3.25,4.25);
      \draw[cell] (3.25,4) rectangle (3.5,4.25);
      \draw[cell] (3.5,4) rectangle (3.75,4.25);
      \draw[cell] (3.75,4) rectangle (4,4.25);
      \draw[cell] (4,4) rectangle (4.25,4.25);
      \draw[cell] (4.25,4) rectangle (4.5,4.25);
      \draw[cell] (4.5,4) rectangle (4.75,4.25);
      \draw[cell] (4.75,4) rectangle (5,4.25);
      \draw[cell] (3,4.25) rectangle (3.25,4.5);
      \draw[cell] (3.25,4.25) rectangle (3.5,4.5);
      \draw[cell] (3.5,4.25) rectangle (3.75,4.5);
      \draw[cell] (3.75,4.25) rectangle (4,4.5);
      \draw[cell] (4,4.25) rectangle (4.25,4.5);
      \draw[cell] (4.25,4.25) rectangle (4.5,4.5);
      \draw[cell] (4.5,4.25) rectangle (4.75,4.5);
      \draw[cell] (4.75,4.25) rectangle (5,4.5);
      \draw[cell] (3,4.5) rectangle (3.25,4.75);
      \draw[cell] (3.25,4.5) rectangle (3.5,4.75);
      \draw[cell] (3.5,4.5) rectangle (3.75,4.75);
      \draw[cell] (3.75,4.5) rectangle (4,4.75);
      \draw[cell] (4,4.5) rectangle (4.25,4.75);
      \draw[cell] (4.25,4.5) rectangle (4.5,4.75);
      \draw[cell] (4.5,4.5) rectangle (4.75,4.75);
      \draw[cell] (4.75,4.5) rectangle (5,4.75);
      \draw[cell] (3,4.75) rectangle (3.25,5);
      \draw[cell] (3.25,4.75) rectangle (3.5,5);
      \draw[cell] (3.5,4.75) rectangle (3.75,5);
      \draw[cell] (3.75,4.75) rectangle (4,5);
      \draw[cell] (4,4.75) rectangle (4.25,5);
      \draw[cell] (4.25,4.75) rectangle (4.5,5);
      \draw[cell] (4.5,4.75) rectangle (4.75,5);
      \draw[cell] (4.75,4.75) rectangle (5,5);
      \draw[dombox] (0,0) rectangle (6,6);
      \draw[ghostline] (1,1) rectangle (3,3);
      \draw[suppline] (2,2) rectangle (5,5); 
      \node[mkD] at (3.375,3.375) {};
      \node[mkD] at (3.625,3.375) {};
      \node[mkD] at (3.875,3.375) {};
      \node[mkD] at (4.125,3.375) {};
      \node[mkD] at (4.375,3.375) {};
      \node[mkD] at (4.625,3.375) {};
      \node[mkD] at (3.375,3.625) {};
      \node[mkD] at (3.625,3.625) {};
      \node[mkD] at (3.875,3.625) {};
      \node[mkD] at (4.125,3.625) {};
      \node[mkD] at (4.375,3.625) {};
      \node[mkD] at (4.625,3.625) {};
      \node[mkD] at (3.375,3.875) {};
      \node[mkD] at (3.625,3.875) {};
      \node[mkD] at (3.875,3.875) {};
      \node[mkD] at (4.125,3.875) {};
      \node[mkD] at (4.375,3.875) {};
      \node[mkD] at (4.625,3.875) {};
      \node[mkD] at (3.375,4.125) {};
      \node[mkD] at (3.625,4.125) {};
      \node[mkD] at (3.875,4.125) {};
      \node[mkD] at (4.125,4.125) {};
      \node[mkD] at (4.375,4.125) {};
      \node[mkD] at (4.625,4.125) {};
      \node[mkD] at (3.375,4.375) {};
      \node[mkD] at (3.625,4.375) {};
      \node[mkD] at (3.875,4.375) {};
      \node[mkD] at (4.125,4.375) {};
      \node[mkD] at (4.375,4.375) {};
      \node[mkD] at (4.625,4.375) {};
      \node[mkD] at (3.375,4.625) {};
      \node[mkD] at (3.625,4.625) {};
      \node[mkD] at (3.875,4.625) {};
      \node[mkD] at (4.125,4.625) {};
      \node[mkD] at (4.375,4.625) {};
      \node[mkD] at (4.625,4.625) {};
      \node[mkC] at (1.75,3.75) {};
      \node[mkC] at (1.75,4.25) {};
      \node[mkC] at (2.25,3.75) {};
      \node[mkC] at (2.25,4.25) {};
      \node[mkC] at (2.75,3.75) {};
      \node[mkC] at (2.75,4.25) {};
      \node[mkC] at (3.25,3.75) {};
      \node[mkC] at (3.25,4.25) {};
      \node[mkC] at (3.75,1.75) {};
      \node[mkC] at (3.75,2.25) {};
      \node[mkC] at (3.75,2.75) {};
      \node[mkC] at (3.75,3.25) {};
      \node[mkC] at (4.25,1.75) {};
      \node[mkC] at (4.25,2.25) {};
      \node[mkC] at (4.25,2.75) {};
      \node[mkC] at (4.25,3.25) {};
      \node[mkAon] at (2.5,2.5) {};
      \node[mkAon] at (3.5,2.5) {};
      \node[mkAon] at (2.5,3.5) {};
      \node[mkAon] at (3.5,3.5) {};
      \node[panellabel] at (3,-0.45) {(b) coarsened configuration $\hat{\mathcal{Q}}$, $\hat{\mathcal{H}}$};
    \end{scope}
 
    \begin{scope}[shift={(0.1,-1.9)}]
        \node[mkAoff] at (0.15,0) {};
        \node[leg] at (0.4,0) {$\mathcal{A}_{\loc}$ (level $\ell$), reactivated};
        \node[mkF] at (6.55,0) {};
        \node[leg] at (6.8,0) {$\mathcal{C}_{\loc}^{\fineonly}$, removed};
        \node[mkC] at (0.15,-0.55) {};
        \node[leg] at (0.4,-0.55) {$\mathcal{C}_{\loc}^{\mathrm{coarse}}$, kept, coefficient corrected};
        \node[mkD] at (6.55,-0.55) {};
        \node[leg] at (6.8,-0.55) {$\mathcal{D}_{\loc}$ (levels $\ge\ell+2$), kept, coefficient corrected};
    \end{scope}
    \end{tikzpicture}
    \caption{Local coarsening of the marked region $\Omega_{\loc}$ (shaded) for
        biquadratic splines with three levels. Basis functions are drawn at their
        Greville abscissae; hollow markers denote deactivated B-splines and solid
        markers active ones. \emph{(a)} The fine configuration $(\QQ,\HH)$. \emph{(b)} The coarse configuration $(\hat\QQ,\hat\HH)$, obtained by coarsening $\Omega_{\loc}$; the ancestors
        $\mathcal{A}_{\loc}$ are reactivated and the functions in
        $\mathcal{C}_{\loc}^{\fineonly}$ disappear, while all remaining
        B-splines keep their status. The dashed square is the support of one
        ancestor: it reaches beyond $\Omega_{\loc}$ and overlaps functions of $\mathcal C_{\loc}^{\mathrm{coarse}}$ and $\mathcal D_{\loc}$, whose coefficients are therefore modified by the coarsening step. B-splines outside these
        four sets are not shown; their coefficients are copied.}
    \label{fig:coarsening_sets}
\end{figure}
 
\paragraph{Block structure of the refinement matrix.} With the notation above, the two bases can be written as the disjoint unions
\begin{equation}\label{eq:hierarchical_basis_disjoint_partitions}
    \hat\HH
    =\mathcal A_{\loc}\,\dot\cup\,
    \mathcal C_{\loc}^{\mathrm{coarse}}\,\dot\cup\,
    \mathcal D_{\loc}\,\dot\cup\,\mathcal H_{\mathrm{rem}},
    \qquad
    \HH
    =\mathcal C_{\loc}^{\fineonly}\,\dot\cup\,
    \mathcal C_{\loc}^{\mathrm{coarse}}\,\dot\cup\,
    \mathcal D_{\loc}\,\dot\cup\,\mathcal H_{\mathrm{rem}}.
\end{equation}
The first component is the genuinely changing part. The ancestors $\mathcal A_{\loc}$ disappear from the fine basis and are replaced by the fine-only functions $\mathcal C_{\loc}^{\fineonly}$. The sets $\mathcal C_{\loc}^{\mathrm{coarse}}$ and $\mathcal D_{\loc}$ remain active, but they are singled out because they appear in the expansion of the ancestors: when a spline is rewritten in the fine basis, their coefficients receive a term from those of $\mathcal A_{\loc}$. In contrast, the functions in $\mathcal H_{\mathrm{rem}}$ neither change their status nor appear in that expansion. This distinction is important for the inverse operation: when coarsening is attempted, only the coefficients attached to the first three local groups have to be recomputed, whereas all coefficients indexed by $\mathcal H_{\mathrm{rem}}$ are copied without any projection.
 
After ordering the rows of the refinement matrix $K$ of~\eqref{eq:refinement_relation} according to the decomposition of $\HH$ and the columns according to the decomposition of $\hat\HH$ in \eqref{eq:hierarchical_basis_disjoint_partitions}, it has the block structure
\begin{equation}\label{eq:global_refinement_detailed_block}
    K=
    \begin{blockarray}{cccc l}
        \mathcal{A}_{\loc} & \mathcal{C}^{\mathrm{coarse}}_{\loc} & \mathcal{D}_{\loc} & \mathcal{H}_{\mathrm{rem}} & \\
        \begin{block}{(cccc)l}
            K_F & 0 & 0 & 0 & \;\mathcal{C}^{\fineonly}_{\loc} \\
            K_C & I & 0 & 0 & \;\mathcal{C}^{\mathrm{coarse}}_{\loc} \\
            K_D & 0 & I & 0 & \;\mathcal{D}_{\loc} \\
            0 & 0 & 0 & I & \;\mathcal{H}_{\mathrm{rem}} \\
        \end{block}
    \end{blockarray},
\end{equation}
where, in terms of the two-level subdivision matrix $C_\ell^{\ell+1}$ of Section~\ref{sec:hierarchical_bsplines},
\begin{equation}\label{eq:coarsening_local_blocks}
    K_F=C_\ell^{\ell+1}(\mathcal C_{\loc}^{\fineonly},\mathcal A_{\loc}),
    \quad
    K_C=C_\ell^{\ell+1}(\mathcal C_{\loc}^{\mathrm{coarse}},\mathcal A_{\loc}),
    \quad
    K_D=M_{\loc}\,C_\ell^{\ell+1}(\mathcal C_{\loc}^{\mathrm{deact}},\mathcal A_{\loc}).
\end{equation}
Equivalently, collecting the first three groups of each basis in a single local block, the refinement matrix has the compact form
\begin{equation}\label{eq:global_refinement_block}
    K=\begin{pmatrix}
        K_{\loc} & 0 \\
        0 & I
    \end{pmatrix},
    \qquad \text{where\,}
    K_{\loc}:=\begin{pmatrix}
        K_F & 0 & 0 \\
        K_C & I & 0 \\
        K_D & 0 & I
    \end{pmatrix}.
\end{equation}
Notice that $K_{\loc}$ contains all the nontrivial entries of $K$, whereas the identity block corresponds to the functions of $\mathcal H_{\mathrm{rem}}$, which have the same coefficient in both representations. Consequently, the refinement operator can be interpreted as a local operator from $\Span\hat\HH_{\loc}$ to $\Span\HH_{\loc}$, where
\begin{equation}\label{eq:local_bases}
    \hat\HH_{\loc}:=\mathcal A_{\loc}\cup\mathcal C_{\loc}^{\mathrm{coarse}}\cup\mathcal D_{\loc},
    \qquad
    \HH_{\loc}:=\mathcal C_{\loc}^{\fineonly}\cup\mathcal C_{\loc}^{\mathrm{coarse}}\cup\mathcal D_{\loc}.
\end{equation}
Throughout we write $n_A:=\#\mathcal A_{\loc}$, $n_F:=\#\mathcal C_{\loc}^{\fineonly}$, $n_C:=\#\mathcal C_{\loc}^{\mathrm{coarse}}$, $n_D:=\#\mathcal D_{\loc}$, $\hat n_{\loc}:=\#\hat\HH_{\loc}=n_A+n_C+n_D$ and $n_{\loc}:=\#\HH_{\loc}=n_F+n_C+n_D$, so that $K_{\loc}\in\RR^{n_{\loc}\times\hat n_{\loc}}$. Since $\mathcal H_{\mathrm{rem}}$ is common to both bases, the partitions~\eqref{eq:hierarchical_basis_disjoint_partitions} also give the dimension of the coarse space in terms of the fine one,
\begin{equation}\label{eq:coarsened_dimension_update}
    \hat n=n-n_F+n_A.
\end{equation}
  In particular, $\hat n\le n$ and $n_A\le n_F$, since $\Span\hat\HH\subseteq\Span\HH$.

Coarsening therefore amounts to realizing a left inverse of $K_{\loc}$. Although $K_{\loc}$ may involve many levels through $\mathcal D_{\loc}$, Section~\ref{S:local_ls_solution} shows that the problem reduces to the fine-only block $K_F$, a submatrix of the subdivision matrix between the two consecutive levels $\ell$ and $\ell+1$.

To illustrate the local sets $\mathcal A_{\loc}$, $\mathcal C_{\loc}^{\fineonly}$, $\mathcal C_{\loc}^{\mathrm{coarse}}$, $\mathcal D_{\loc}$ and the block structure~\eqref{eq:global_refinement_detailed_block}, let us consider the following example.
 
\begin{example}
    We consider biquadratic splines ($p=2$). Figure~\ref{fig:meshes} displays the initial mesh on the left and the resulting mesh on the right, where only a single cell at level $0$ has been refined. This single refinement deactivates exactly one B-spline $\beta_0$ of level $0$, which is therefore the only ancestor, $\mathcal A_{\loc}=\{\beta_0\}$, and activates four new functions at level $1$ ($\mathcal{C}_{\text{loc}}^{\fineonly}$). In this specific scenario, two children of $\beta_0$ were already active in the initial basis ($\mathcal{C}_{\text{loc}}^{\text{coarse}}$). Furthermore, the set $\mathcal{D}_{\text{loc}}$, consisting of the active descendants of the deactivated children of $\beta_0$, contains $8$ functions at level $2$ and $48$ at level $3$, that is, $56$ multi-level descendants in total. The support of $\beta_0$ is $[\tfrac14,1]\times[0,\tfrac14]$, so that it extends well beyond the refined cell, and its expansion in the fine basis also involves functions that do not overlap that cell, namely those of $\mathcal C_{\loc}^{\mathrm{coarse}}$ and $\mathcal D_{\loc}$. Overall the mesh grows from $103$ to $106$ active cells and the space dimension from $\hat n=123$ to $n=126$, in agreement with \eqref{eq:coarsened_dimension_update}, since $126-4+1=123$.
 
    \begin{figure}[htbp]
        \centering
        \includegraphics[width=0.49\textwidth]{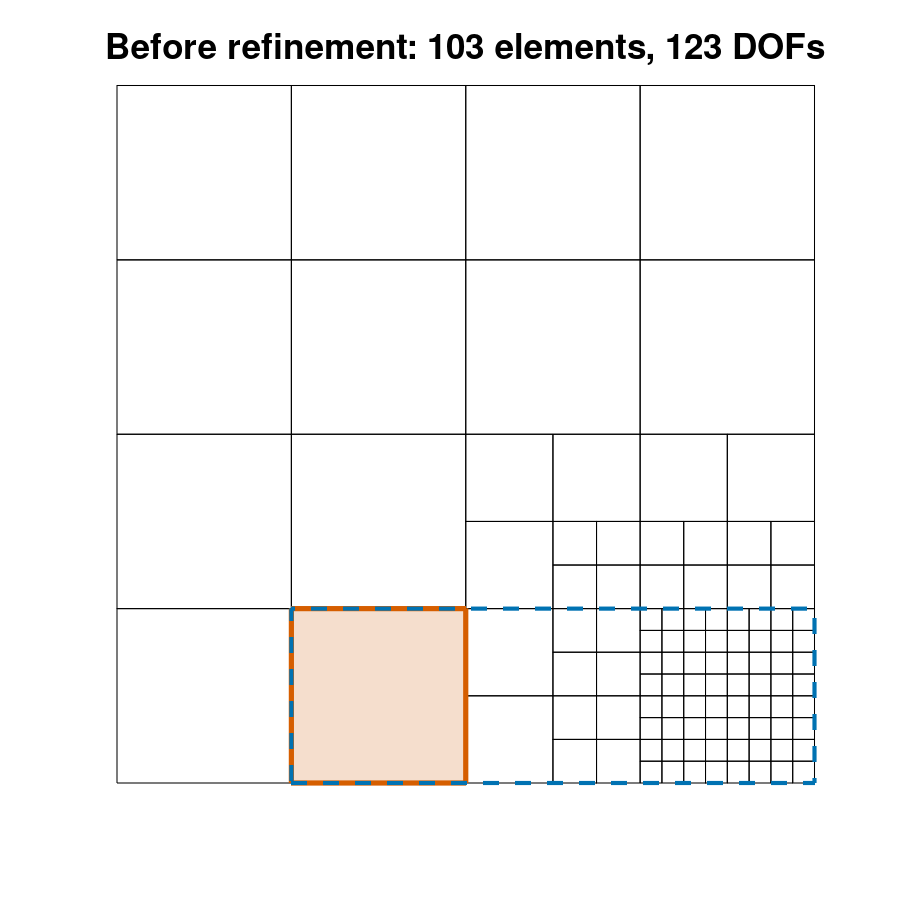}
        \hfill
        \includegraphics[width=0.49\textwidth]{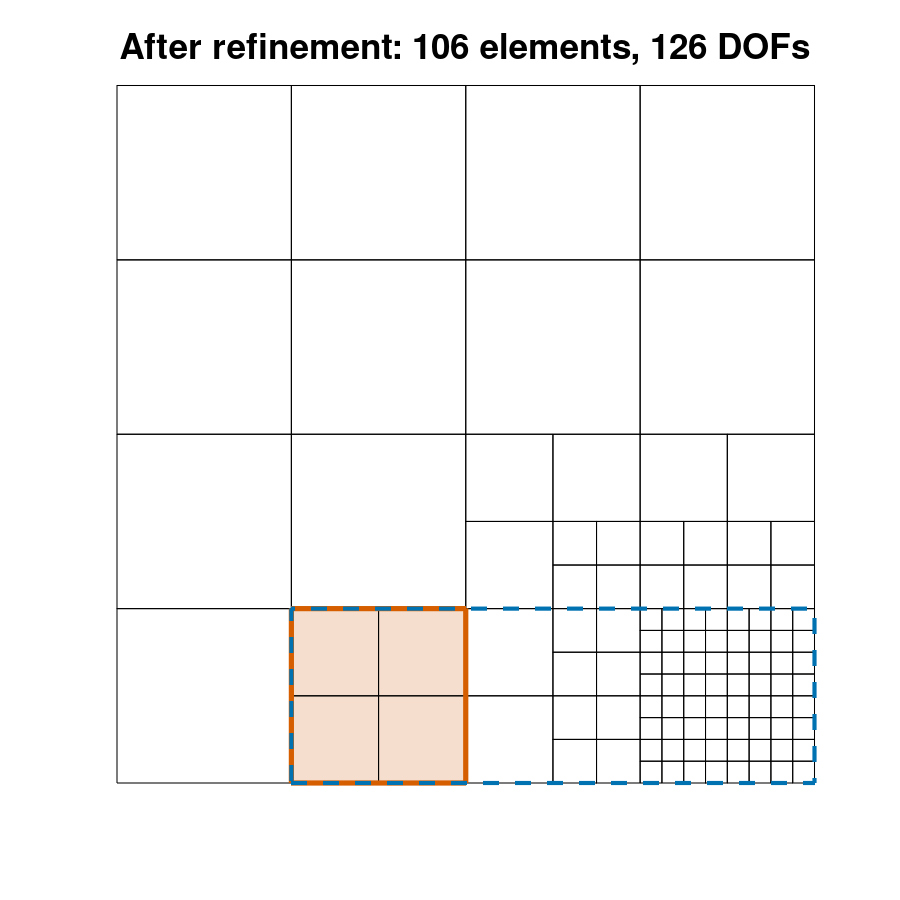}
        \caption{Left: Initial hierarchical mesh. Right: Refined mesh after subdivision of a single cell at level 0. The shaded cell is the marked region $\Omega_{\loc}$, and the dashed region is the support of the B-spline $\beta_0$, which is active in the initial mesh and deactivated by the refinement.}
        \label{fig:meshes}
    \end{figure}

    The sparsity pattern of the corresponding refinement matrix is depicted in Figure~\ref{fig:sparsity}, which isolates the column corresponding to $\beta_0$ and highlights the fine basis functions involved in its expansion. The right panel displays the entries of that single column, which are the coefficients of the expansion of $\beta_0$ in the fine basis: four large entries correspond to the new level-$1$ functions of $\mathcal C_{\loc}^{\fineonly}$, two to the pre-existing ones of $\mathcal C_{\loc}^{\mathrm{coarse}}$, and the remaining $56$, those of $\mathcal D_{\loc}$, are spread over levels $2$ and $3$, decaying rapidly with the level. In the coarsening step of Section~\ref{S:local_coarsening_problem} below, the coefficient of $\beta_0$ is recovered from the four coefficients of $\mathcal C_{\loc}^{\fineonly}$ (here $K_F$ is a $4\times1$ matrix), and the term that the expansion of $\beta_0$ adds to the remaining $2+56$ coefficients is subtracted.
 
    \begin{figure}[htbp]
        \centering
        \includegraphics[width=0.49\textwidth]{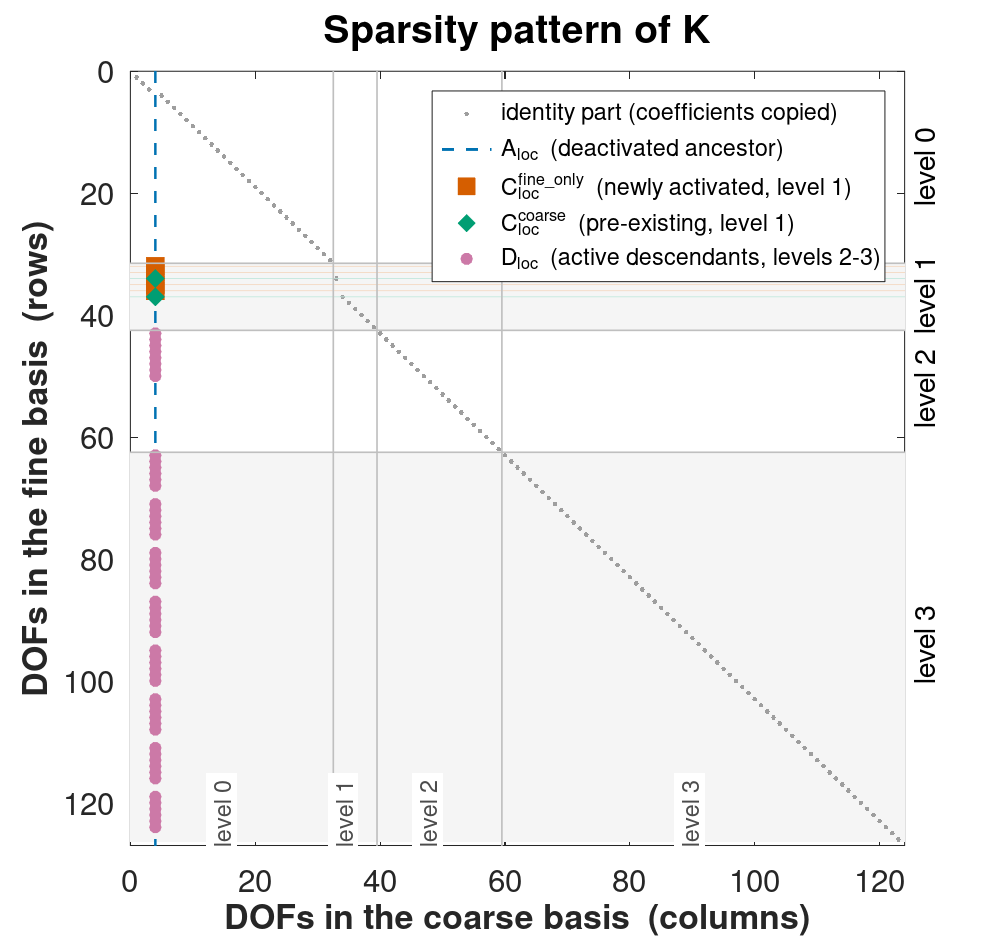}
        \hfill
        \includegraphics[width=0.49\textwidth]{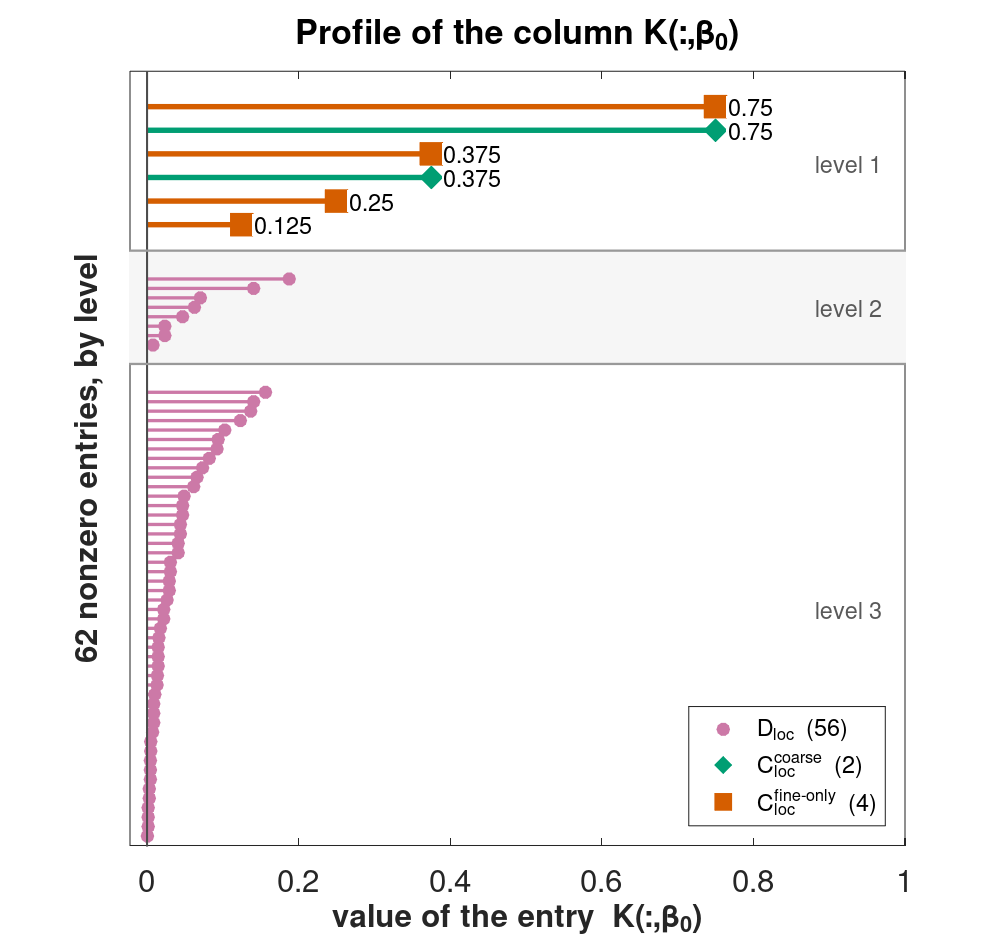}
        \caption{Left: sparsity pattern of the refinement matrix, with the degrees of freedom grouped by level and the column of $\beta_0$ highlighted; away from that column the matrix reduces to an identity block, which is broken into level blocks by the different level ordering of the two bases. Right: the $62$ nonzero entries of that column, grouped by level and sorted by decreasing magnitude within each group. Of the six level-$1$ entries, four correspond to $\mathcal C_{\loc}^{\fineonly}$ and two to $\mathcal C_{\loc}^{\mathrm{coarse}}$; the $56$ entries of levels $2$ and $3$ correspond to $\mathcal D_{\loc}$. The role of each group in the coarsening step is explained in Section~\ref{S:local_ls_solution}. Note that the entries are subdivision weights and do not add up to one.}
        \label{fig:sparsity}
    \end{figure}
\end{example}

\subsection{The local coarsening problem and an error bound}\label{S:local_coarsening_problem}\label{S:local_ls_solution}
 
We now reverse the point of view. We start from a hierarchical mesh $\QQ$, its hierarchical basis $\HH$, and a subset of active fine cells $\QQ_{\loc} \subset \GG_{\ell+1}\cap\QQ$ marked for coarsening. The aim is to decide whether a spline of $\Span\HH$ can be represented, up to a prescribed tolerance, in the hierarchical basis of the coarser mesh in which these cells are replaced by their parents.
 
Since a cell of level $\ell$ can only be reactivated if all its $2^d$ children are removed, we require the marked cells to form complete sibling groups,
\begin{equation}\label{eq:sibling_groups_coarsening}
    \QQ_{\loc} = \children{\hat\QQ_{\loc}},
\end{equation}
where $\hat\QQ_{\loc}:=\parents{\QQ_{\loc}}\subset\GG_\ell$; marked cells that do not complete a sibling group are discarded. 
 
The coarse hierarchical mesh $\hat{\QQ}$ is defined by the hierarchy of subdomains $\hat{\boldsymbol\Omega}_L=\{\hat{\Omega}_{k}\}_{k=0}^L$ given by
\begin{equation}\label{eq:coarse_hierarchy}
    \hat\Omega_{\ell+1}:=\overline{\Omega_{\ell+1}\setminus\Omega_{\loc}},
    \qquad
    \hat\Omega_k:=\Omega_k\ \ (k\neq\ell+1),
\end{equation}
where $\Omega_{\loc}$ is the region~\eqref{eq:local_region}. These subdomains indeed form a hierarchy: $\hat\Omega_{\ell+1}$ is a union of cells of $\GG_\ell$, because so are $\Omega_{\ell+1}$ and $\Omega_{\loc}$, and it contains $\Omega_{\ell+2}$, because the marked cells are active. Moreover, refining the cells in $\hat\QQ_{\loc}$ recovers $\{\Omega_k\}_{k=0}^L$, cf.~\eqref{eq:refined_hierarchy}. Hence the coarse basis $\hat\HH$ and the fine basis $\HH$ are exactly related as in Section~\ref{S:local_block_structure}: the local sets \eqref{eq:def_Aloc}--\eqref{eq:def_Dloc} are computed from the fine configuration alone, without constructing $\hat\HH$ beforehand, and the local block $K_{\loc}$ is given by \eqref{eq:coarsening_local_blocks}--\eqref{eq:global_refinement_block}.

Given $s=\sum_{\beta\in\HH}c_\beta\,\beta\in\Span\HH$, let $c_F$, $c_C$ and $c_D$ denote the restrictions of its coefficient vector $c$ to $\mathcal C_{\loc}^{\fineonly}$, $\mathcal C_{\loc}^{\mathrm{coarse}}$ and $\mathcal D_{\loc}$, so that $c_{\loc}:=(c_F,c_C,c_D)^T$ is its restriction to $\HH_{\loc}$, ordered as in~\eqref{eq:local_bases}. The coefficient vector $\hat c\in\RR^{\hat n}$ of the coarse spline coincides with $c$ on $\mathcal H_{\mathrm{rem}}$, and its restriction $\hat c_{\loc}=(\hat c_A,\hat c_C,\hat c_D)^T$ to $\hat\HH_{\loc}$, with blocks indexed by $\mathcal A_{\loc}$, $\mathcal C_{\loc}^{\mathrm{coarse}}$ and $\mathcal D_{\loc}$, is defined as a solution of the \emph{local discrete least-squares problem}
\begin{equation}\label{eq:full_ls_problem}
    \hat c_{\loc}=\argmin_{v\in\RR^{\hat n_{\loc}}}\ \|c_{\loc}-K_{\loc}v\|_2.
\end{equation}
 In block form, \eqref{eq:full_ls_problem} reads
\begin{equation*}
    \begin{pmatrix} c_F \\ c_C \\ c_D \end{pmatrix} \approx \begin{pmatrix} K_F & 0 & 0 \\ K_C & I & 0 \\ K_D & 0 & I \end{pmatrix} \begin{pmatrix} \hat{c}_A \\ \hat{c}_C \\ \hat{c}_D \end{pmatrix}.
\end{equation*}
Two features of \eqref{eq:full_ls_problem} deserve a comment. First, the minimization is performed in the Euclidean norm of the \emph{coefficient} vectors: no mass matrix and no quadrature are involved, and the result is in general not the $L^2(\Omega)$-projection onto $\Span\hat\HH$. Second, the procedure only modifies the spline near the region where the mesh changes: as shown below, its result belongs to the class
\begin{equation}\label{eq:admissible_class}
    \hat{\mathcal V}:=\big\{\hat v\in\Span\hat\HH:\ \text{the coefficients of } \hat v \text{ and } s \text{ in } \HH \text{ coincide on } \mathcal{F}^A\setminus\mathcal C_{\loc}^{\fineonly}\big\},
\end{equation}
and \eqref{eq:full_ls_problem} selects the element of $\hat{\mathcal V}$ whose fine coefficients on $\mathcal C_{\loc}^{\fineonly}$ are closest to those of $s$. Since the elements of $\hat{\mathcal V}$ keep all the other fine coefficients of $s$, the coarse spline is in general less accurate in the $L^2$-norm than the $L^2(\Omega)$-projection, which may modify every coefficient; this is the price of locality, quantified in Section~\ref{sec:numerical_experiments}. In exchange, the error admits the computable bound of Proposition~\ref{prop:linf_bound} below. We now show that this element is unique, that it is determined by a least-squares problem for the fine-only block $K_F$ alone, and that its residual measures exactly the information lost by the coarsening step.
 
We first establish a key algebraic property of the fine-only block $K_F$.
 
\begin{lemma}[Full column rank of $K_F$]\label{lemma:kf_rank}
    The matrix $K_F\in\RR^{n_F\times n_A}$ of~\eqref{eq:coarsening_local_blocks} has full column rank.
\end{lemma}

\begin{proof}
    We first show that $K_{\loc}$ has full column rank. By~\eqref{eq:refinement_relation}, for every $\hat v\in\RR^{\hat n}$, $K\hat v$ is the coefficient vector in $\HH$ of the spline $\sum_{\hat\beta\in\hat\HH}\hat v_{\hat\beta}\hat\beta$. If $K\hat v=\mathbf 0$, this spline vanishes, and $\hat v=\mathbf 0$ because $\hat\HH$ is a basis. Hence $K$ has full column rank, and so does its diagonal block $K_{\loc}$ in~\eqref{eq:global_refinement_block}.

    Now let $x_A\in\RR^{n_A}$ satisfy $K_Fx_A=\mathbf 0$. The vector $x=(x_A,-K_Cx_A,-K_Dx_A)^T$ satisfies $K_{\loc}x=\mathbf 0$ by the block form~\eqref{eq:global_refinement_block}, hence $x=\mathbf 0$ and in particular $x_A=\mathbf 0$.
\end{proof}

\begin{theorem}[Solution and residual of the local problem]\label{thm:exact_residual}
    The local least-squares problem \eqref{eq:full_ls_problem} has a unique solution $\hat{c}_{\loc} = (\hat{c}_A, \hat{c}_C, \hat{c}_D)^T$, given by $\hat{c}_A = K_F^\dagger c_F$, where $K_F^\dagger=(K_F^TK_F)^{-1}K_F^T$ is the Moore--Penrose pseudoinverse of $K_F$, and
    \begin{align}
        \hat{c}_C & = c_C - K_C \, \hat{c}_A, \label{eq:exact_cc} \\
        \hat{c}_D & = c_D - K_D \, \hat{c}_A. \label{eq:exact_cd}
    \end{align}
    Moreover, the residual $r_{\loc} := c_{\loc} - K_{\loc} \hat{c}_{\loc}$ vanishes on $\mathcal C_{\loc}^{\mathrm{coarse}}$ and $\mathcal D_{\loc}$, and coincides on $\mathcal C_{\loc}^{\fineonly}$ with the residual of the reduced system $K_F\hat c_A\approx c_F$; that is,
    \begin{equation}\label{eq:collapsed_residual}
        r_{\loc} = \begin{pmatrix} c_F - K_F \, \hat{c}_A \\ 0 \\ 0 \end{pmatrix}.
    \end{equation}
\end{theorem}

\begin{proof}
    For $v=(v_A,v_C,v_D)^T\in\RR^{\hat n_{\loc}}$, the definition of $K_{\loc}$ in~\eqref{eq:global_refinement_block} gives
    \[
        \|c_{\loc}-K_{\loc}v\|_2^2=\|c_F-K_Fv_A\|_2^2+\|c_C-K_Cv_A-v_C\|_2^2+\|c_D-K_Dv_A-v_D\|_2^2.
    \]
    For every $v_A$, the last two terms are nonnegative and vanish if and only if $v_C=c_C-K_Cv_A$ and $v_D=c_D-K_Dv_A$. Hence $v$ minimizes $\|c_{\loc}-K_{\loc}v\|_2$ if and only if $v_A$ minimizes $\|c_F-K_Fv_A\|_2$ and $v_C$, $v_D$ are given by these formulas. Since $K_F$ has full column rank by Lemma~\ref{lemma:kf_rank}, the normal equations $K_F^TK_Fv_A=K_F^Tc_F$ have the unique solution $v_A=K_F^\dagger c_F$. This proves uniqueness and~\eqref{eq:exact_cc}--\eqref{eq:exact_cd}; for this choice the last two blocks of $c_{\loc}-K_{\loc}v$ vanish, which is~\eqref{eq:collapsed_residual}.
\end{proof}

We call~\eqref{eq:exact_cc}--\eqref{eq:exact_cd} the \emph{counter-corrections}: they subtract from $c_C$ and $c_D$ the terms $K_C\hat c_A$ and $K_D\hat c_A$ that the expansion of the ancestors in the fine basis adds to them.

Let $\Pi:\Span\HH\to\Span\hat\HH$ denote the resulting coarsening operator, which maps each $s\in\Span\HH$ to the spline with coefficient vector $\hat c$ in $\hat\HH$; it is linear, since $\hat c$ depends linearly on $c$. The following consequences of Theorem~\ref{thm:exact_residual} follow at once.

\begin{corollary}[Left inverse and projection]\label{cor:exactness}\label{cor:compression}
    The map $c\mapsto\hat c$ is a left inverse of $K$, and $c_{\loc}\mapsto\hat c_{\loc}$ is a left inverse of $K_{\loc}$. Consequently, if $s\in\Span\hat\HH$, then $\Pi s=s$ and the residual $r_{\loc}$ of Theorem~\ref{thm:exact_residual} vanishes; in particular, $\Pi$ is a projection onto $\Span\hat\HH$.
\end{corollary}

\begin{proof}
    By Theorem~\ref{thm:exact_residual}, $\hat c_{\loc}=Xc_{\loc}$, and since $\hat c$ is copied from $c$ on $\mathcal H_{\mathrm{rem}}$, $\hat c=Yc$, where
    \[
        X:=\begin{pmatrix} K_F^\dagger & 0 & 0 \\ -K_C K_F^\dagger & I & 0 \\ -K_D K_F^\dagger & 0 & I \end{pmatrix},
        \qquad
        Y:=\begin{pmatrix} X & 0 \\ 0 & I \end{pmatrix},
    \]
    with blocks ordered as in~\eqref{eq:global_refinement_block}. Since $K_F^\dagger K_F=I$, block multiplication with~\eqref{eq:global_refinement_block} gives $XK_{\loc}=I$, and hence $YK=I$. Now let $s\in\Span\hat\HH$, with coefficient vector $\hat v$ in $\hat\HH$. By~\eqref{eq:refinement_relation}, $c=K\hat v$, so that $\hat c=YK\hat v=\hat v$, that is, $\Pi s=s$. Therefore $c=K\hat c$, and~\eqref{eq:global_refinement_block} gives $c_{\loc}=K_{\loc}\hat c_{\loc}$, whence $r_{\loc}=\mathbf 0$.
\end{proof}

For a general $s\in\Span\HH$, the residual provides a guaranteed error bound, computable at negligible cost and before the mesh is modified: since the nonzero entries of $r_{\loc}$ are indexed by $\mathcal C_{\loc}^{\fineonly}$, a set of B-splines of a single level, $\|r_{\loc}\|_\infty$ bounds the $L^\infty$-error of the coarsening step with no unknown constants.

\begin{proposition}[Constant-free $L^\infty$ error bound]\label{prop:linf_bound}
    Let $s\in\Span\HH$ and $\hat s:=\Pi s$. Then
    \begin{equation}\label{eq:single_level_inf_bound}
        s-\hat s=\sum_{\beta\in\mathcal C_{\loc}^{\fineonly}}(r_{\loc})_\beta\,\beta
        \qquad\text{and}\qquad
        \|s-\hat s\|_{L^\infty(\Omega)}\le\|r_{\loc}\|_\infty.
    \end{equation}
    In particular, $s-\hat s$ vanishes identically outside $\bigcup\{\supp\beta:\beta\in\mathcal C_{\loc}^{\fineonly}\}$.
\end{proposition}
\begin{proof}
    By~\eqref{eq:refinement_relation} and~\eqref{eq:global_refinement_block}, the coefficients of $s-\hat s$ in $\HH$ form the vector $c-K\hat c$, which coincides with $r_{\loc}$ on $\HH_{\loc}$ and vanishes on $\mathcal H_{\mathrm{rem}}$, where $\hat c$ is copied from $c$. Together with~\eqref{eq:collapsed_residual}, this proves the identity in~\eqref{eq:single_level_inf_bound}, and hence the last statement. The functions of $\mathcal C_{\loc}^{\fineonly}$ belong to $\mathcal B_{\ell+1}$, whose B-splines are non-negative and form a partition of unity; therefore $\sum_{\beta\in\mathcal C_{\loc}^{\fineonly}}\beta(x)\le1$ and
    \[
        |s(x)-\hat s(x)|\le\max_{\beta}|(r_{\loc})_\beta|\sum_{\beta\in\mathcal C_{\loc}^{\fineonly}}\beta(x)\le\|r_{\loc}\|_\infty
    \]
    for every $x\in\Omega$.
\end{proof}
 
\begin{remark}[Computational cost]\label{rem:complexity}
    The whole coarsening step reduces to the thin overdetermined system $K_F\hat c_A\approx c_F$, of size $n_F\times n_A$ with $n_A\le n_F$, followed by the two matrix--vector counter-corrections \eqref{eq:exact_cc}--\eqref{eq:exact_cd}. Solving it by a QR factorization costs $\mathcal O(n_F n_A^2)$ operations. The entries of $K_F$ are read off the two-level subdivision matrix $C_\ell^{\ell+1}$, which is available in standard hierarchical spline libraries as a Kronecker product of univariate subdivision matrices of the $d$ parametric directions, and each of its columns has at most $(p+2)^d$ nonzero entries---the number of children of a single B-spline under dyadic refinement---regardless of the size of the coarsening region. When $\hat\QQ_{\loc}$ consists of a single cell, $n_A\le(p+1)^d$ and $n_F\le(p+2)^d$; both bounds are attained when all the level-$\ell$ B-splines that do not vanish on the cell are deactivated, since then all the level-$(\ell+1)$ B-splines that do not vanish on it are active, and in general $K_F$ is a submatrix of the corresponding $(p+2)^d\times(p+1)^d$ block of $C_\ell^{\ell+1}$. The counter-corrections cost $\mathcal O(\operatorname{nnz}(K_C)+\operatorname{nnz}(K_D))$ operations, plus the assembly of $K_D$ through $M_{\loc}$; only this part grows with the number of finer levels present in the supports of the ancestors, through $n_D$. No mass matrix or global system is ever assembled, and no quadrature is needed.
\end{remark}
 
 \begin{remark}[Coarsening without ancestors]\label{rem:no_ancestors}
    A B-spline $\gamma\in\mathcal B_\ell$ is deactivated only if its support, $p+1$ cells of level $\ell$ wide in each direction, fits inside $\Omega_{\ell+1}$, whereas $\beta\in\mathcal B_{\ell+1}$, whose support is half as wide, is active as soon as $\supp\beta\subset\Omega_{\ell+1}$ and $\supp\beta\not\subset\Omega_{\ell+2}$. A refined strip whose width is at least $(p+1)/2$ but less than $p+1$ cells of level $\ell$---two cells for biquadratic splines, the typical outcome of tracking a moving front---therefore carries active functions of level $\ell+1$ and no deactivated function of level $\ell$, so that a legitimate coarsening step may face $\mathcal A_{\loc}=\varnothing$ while $\mathcal C_{\loc}^{\fineonly}\neq\varnothing$.

    Since Definition~\eqref{eq:def_fineonly} does not require the fine-only functions to be children of the ancestors, this case needs no special treatment. A fine-only function that is not a child of any ancestor has a zero row in $K_F$, so by~\eqref{eq:collapsed_residual} its entry in $r_{\loc}$ equals its coefficient, which thus enters the bound of Proposition~\ref{prop:linf_bound} in full. In particular, when $\mathcal A_{\loc}=\varnothing$ one has $K_F\in\RR^{n_F\times0}$ and the vector $\hat c_A$ is void: nothing can be recovered, and the error is the whole part of $s$ carried by the removed functions, $s-\hat s=\sum_{\beta\in\mathcal C_{\loc}^{\fineonly}}c_\beta\,\beta$.
\end{remark}

\begin{remark}[Other left inverses]\label{prop:canonical_left_inverse}\label{rem:left_inverse_family}\label{rem:alternative_left_inverses}
    Replacing $K_F^\dagger$ by any left inverse $K_F^L$ of $K_F$ in the matrix $X$ of the proof of Corollary~\ref{cor:exactness} yields, by the same block multiplication, a matrix $X_L$ with $X_LK_{\loc}=I$. If $\hat c_{\loc}$ is defined as $X_Lc_{\loc}$, Corollary~\ref{cor:exactness} and Proposition~\ref{prop:linf_bound} remain valid, since the proof of the corollary only uses $K_F^LK_F=I$ and the residual keeps the structure~\eqref{eq:collapsed_residual}, with first block $(I-K_FK_F^L)c_F$. What changes is the norm of the residual: $K_F^\dagger$ minimizes its Euclidean norm, but not necessarily its maximum norm, which is the one that enters the bound; see Remark~\ref{rem:l2_vs_linf}. Cheaper choices of $K_F^L$, in which the row of each $\beta\in\mathcal A_{\loc}$ is an explicit local functional $\lambda_\beta$, so that $(\hat c_A)_\beta=\lambda_\beta(c_F)$ and no least-squares problem has to be solved, could be derived from the quasi-interpolants for hierarchical spaces of~\cite{SM2016,Speleers2017} or from the local left inverses of tensor-product subdivision matrices of~\cite{AFG25}. Neither construction applies verbatim, since the functionals involved may use coefficients other than those in $c_F$; adapting them is left for future work.
\end{remark}

\begin{remark}[Euclidean versus maximum norm]\label{rem:l2_vs_linf}
   Problem \eqref{eq:full_ls_problem} minimizes the Euclidean norm of the residual, whereas the guaranteed bound \eqref{eq:single_level_inf_bound} is expressed in the maximum norm, the two being equivalent only up to a factor $\sqrt{n_F}$. The bound holds for the computed $\hat c_A$, but this $\hat c_A$ does not in general minimize $\|c_F-K_F\hat c_A\|_\infty$; hence, in the adaptive algorithm of Section~\ref{sec:adaptive_coarsening}, a step may be rejected for a prescribed tolerance that another choice of $\hat c_A$ would meet. The mismatch thus costs compression, not accuracy. Solving the local Chebyshev problem $\min_{v_A}\|c_F-K_Fv_A\|_\infty$ would close the gap, but it requires a linear program at every evaluation of the bound, instead of the closed-form solution of Theorem~\ref{thm:exact_residual}, and yields in general a nonlinear coarsening operator. We therefore keep the least-squares formulation. Moreover, a rejected step is not entirely lost, since the rollback mechanism of that algorithm retries it on a smaller candidate region. 
\end{remark}

\subsection{An algebraic algorithm for multi-level coarsening}
\label{sec:algorithmic_implementation}
 
The procedure is organized in three routines. The driver, Algorithm~\ref{alg:coarsen_multilevel}, sweeps the levels from the finest transition $\ell=L-2$ down to $\ell=0$; at each level with marked cells, it calls Algorithm~\ref{alg:estimate_error}, which classifies the local functions and evaluates the residual, and then Algorithm~\ref{alg:coarsen_at_level}, which updates the mesh, the basis and the coefficients. We assume that the subdivision matrices $C_\ell^{\ell+1}$ and the transfer matrices $M_\ell$ are available; in GeoPDEs~\cite{geopdes,GV18}, the former are stored through their univariate factors, and the latter are assembled from them by the recursion \eqref{E:multilevel transfer matrices}. To keep the interfaces compact, we collect the local function sets identified by Algorithm~\ref{alg:estimate_error} into the single structure
\begin{equation}\label{eq:local_structure}
    \mathcal L_{\loc}:=\big(\mathcal A_{\loc},\ \mathcal C_{\loc}^{\fineonly}\big),
\end{equation}
which is passed unchanged to Algorithm~\ref{alg:coarsen_at_level}.

Since the marked cells are active cells of the initial mesh forming complete sibling groups, the marked sets of different levels are disjoint and are not affected by the steps performed at the other levels. The levels could therefore be processed in any order, although the resulting spline depends on that order in general. We choose the descending order because it is the natural one for the adaptive algorithm of Section~\ref{sec:adaptive_coarsening}, in which the candidate cells of level $\ell$ are determined on the current mesh, and coarsening the finer levels first enlarges the set of candidates.
 
\begin{algorithm}[h!tb]
    \setstretch{1.25}
    \caption{\bf coarsen\_multilevel}\label{alg:coarsen_multilevel}
    \begin{algorithmic}[1]
        \Statex \textbf{Input:} Hierarchical mesh $\mathcal{Q}$, hierarchical basis $\mathcal{H}$, coefficient vector $c$, depth $L$ of the hierarchy, and marked cells $\mathcal{Q}^{(\ell+1)}_{\loc}\subset\mathcal{Q}\cap\GG_{\ell+1}$ satisfying~\eqref{eq:sibling_groups_coarsening} for $\ell=0,\dots,L-2$ 
        \Statex \textbf{Output:} Coarsened mesh $\mathcal{Q}$, coarsened basis $\mathcal{H}$, updated coefficient vector $c$, and accumulated bound $e_{\infty}$ for the $L^\infty$-error
 
        \State $e_{\infty} \gets 0$ \Comment{Initialize the accumulated bound}
 
        \For{$\ell = L-2$ \textbf{down to} $0$}
        \State $\mathcal{Q}_{\loc} \gets \mathcal{Q}^{(\ell+1)}_{\loc}$
        \If{$\mathcal{Q}_{\loc} \neq \varnothing$}
        \State $(e_{\ell}, \mathcal{L}_{\loc}, \hat{c}_A, r_{\loc}) \gets \text{\bf local\_coarsening\_residual}(\mathcal{Q}, \mathcal{H}, c, \mathcal{Q}_{\loc}, \ell)$ 
 
        \State $e_{\infty} \gets e_{\infty} + e_{\ell}$ \Comment{Accumulate the bounds across levels}
 
        \State $(\mathcal{Q}, \mathcal{H}, c) \gets \text{\bf coarsen\_at\_level}(\mathcal{Q}, \mathcal{H}, c, \mathcal{Q}_{\loc}, \ell, \mathcal{L}_{\loc}, \hat{c}_A)$
        \EndIf
        \EndFor
        \State \textbf{return} $\mathcal{Q}$, $\mathcal{H}$, $c$, $e_{\infty}$
    \end{algorithmic}
\end{algorithm}

 Algorithm~\ref{alg:estimate_error} computes, before any change to the mesh, the local structure $\mathcal L_{\loc}$ of~\eqref{eq:local_structure}, the matrix $K_F$ of~\eqref{eq:coarsening_local_blocks}, the ancestor coefficients $\hat c_A$ and the residual $r_{\loc}$, of which only the nonzero block in~\eqref{eq:collapsed_residual} is computed. By Proposition~\ref{prop:linf_bound}, the returned value $e_{\loc}=\|r_{\loc}\|_\infty$ is a guaranteed bound for the $L^\infty$-error of the step, not an estimate.

\begin{algorithm}[h!tb]
    \setstretch{1.25}
    \caption{\bf local\_coarsening\_residual}\label{alg:estimate_error}
    \begin{algorithmic}[1]
        \Statex \textbf{Input:} Hierarchical mesh $\mathcal{Q}$, hierarchical basis $\mathcal{H}$, coefficient vector $c$, marked cells $\mathcal{Q}_{\loc} \subset \mathcal{Q}\cap\GG_{\ell+1}$ forming complete sibling groups, cf.~\eqref{eq:sibling_groups_coarsening}, and level $\ell$ 
        \Statex \textbf{Output:} $L^\infty$-error bound $e_{\loc}$, local structure $\mathcal{L}_{\loc}$ as in \eqref{eq:local_structure}, coarse ancestor coefficients $\hat{c}_A$, and residual vector $r_{\loc}$
 
        \State $\hat{\mathcal{Q}}_{\loc} \gets \parents{\mathcal{Q}_{\loc}}$, \quad $\Omega_{\loc} \gets \bigcup_{Q \in \hat{\mathcal{Q}}_{\loc}} Q$
        \State $\mathcal{A}_{\loc} \gets \{\beta \in \mathcal{F}_\ell^D : |\supp\beta \cap \Omega_{\loc}| > 0\}$ \Comment{Reactivated coarse ancestors}
        \State $\mathcal{C}_{\loc}^{\fineonly} \gets \{\beta \in \mathcal{F}_{\ell+1}^A: |\supp\beta \cap \Omega_{\loc}| > 0\}$ \Comment{Functions becoming passive}
        \State $K_F \gets C_\ell^{\ell+1}(\mathcal{C}_{\loc}^{\fineonly}, \, \mathcal{A}_{\loc})$
        \State $\hat{c}_A \gets K_F^\dagger \, c(\mathcal{C}_{\loc}^{\fineonly})$ \Comment{Solve the local two-level least-squares problem}
        \State $r_{\loc} \gets c(\mathcal{C}_{\loc}^{\fineonly}) - K_F \, \hat{c}_A$ \Comment{Nonzero block of the residual, cf.~\eqref{eq:collapsed_residual}}
        \State $e_{\loc} \gets \| r_{\loc} \|_\infty$ 
 
        \State \textbf{return} $e_{\loc}$, $\mathcal{L}_{\loc}$, $\hat{c}_A$, $r_{\loc}$
    \end{algorithmic}
\end{algorithm}
 
Algorithm \ref{alg:coarsen_at_level} performs the step. It determines, among the children of the ancestors, the sets $\mathcal C_{\loc}^{\mathrm{coarse}}$ and $\mathcal C_{\loc}^{\mathrm{deact}}$, and the active descendants $\mathcal D_{\loc}$, applies the counter-corrections \eqref{eq:exact_cc}--\eqref{eq:exact_cd} to $\mathcal C_{\loc}^{\mathrm{coarse}}$ and $\mathcal D_{\loc}$, updates the mesh and the basis, and assembles the coarse coefficient vector. The update of the mesh and the basis is performed by \textbf{coarsen\_mesh\_and\_space}, which replaces the marked cells by their parents, as in~\eqref{eq:coarse_hierarchy}, and updates the hierarchical basis accordingly, following the coarsening algorithms of~\cite{GV18} implemented in GeoPDEs~\cite{geopdes}; it does not involve the coefficients. Without the counter-corrections, the term that the expansion of the ancestors adds to the coefficients of $\mathcal C_{\loc}^{\mathrm{coarse}}$ and $\mathcal D_{\loc}$ would change the spline on the supports of these functions, which do not overlap $\Omega_{\loc}$ (see the dashed support of an ancestor in Figure~\ref{fig:coarsening_sets}); with them, the spline is unchanged outside the supports of the removed functions (cf.~Proposition~\ref{prop:linf_bound}). 

\begin{algorithm}[h!tb]
    \setstretch{1.25}
    \caption{\bf coarsen\_at\_level}\label{alg:coarsen_at_level}
    \begin{algorithmic}[1]
        \Statex \textbf{Input:} Hierarchical mesh $\mathcal{Q}$, hierarchical basis $\mathcal{H}$, coefficient vector $c$, marked cells $\mathcal{Q}_{\loc} \subset \mathcal{Q}\cap\GG_{\ell+1}$ forming complete sibling groups, cf.~\eqref{eq:sibling_groups_coarsening}, level $\ell$, local structure $\mathcal{L}_{\loc}$, and ancestor coefficients $\hat{c}_A$.
        \Statex \textbf{Output:} Coarsened mesh $\hat{\mathcal{Q}}$, coarsened basis $\hat{\mathcal{H}}$, and updated coefficient vector $\hat{c}$.
 
        \Comment{Step 1: Classification of the children of the ancestors}
        \State $\mathcal{C}_{\loc} \gets \children{\mathcal{A}_{\loc}}$
        \State $\mathcal{C}_{\loc}^{\mathrm{fine}} \gets \mathcal{C}_{\loc} \cap \mathcal{F}_{\ell+1}^A$
        \State $\mathcal{C}_{\loc}^{\mathrm{coarse}} \gets \mathcal{C}_{\loc}^{\mathrm{fine}} \setminus \mathcal{C}_{\loc}^{\fineonly}$
        \State $\mathcal{C}_{\loc}^{\mathrm{deact}} \gets \mathcal{C}_{\loc} \cap \mathcal{F}_{\ell+1}^D$
 
        \Comment{Step 2: Multi-level transfer to the active descendants}
        \State $\mathcal{D}_{\loc} \gets \{\gamma \in \mathcal{F}_{\ell+2}^A \cup \cdots \cup \mathcal{F}_{L-1}^A : M_{\ell+1}(\{\gamma\}, \mathcal{C}_{\loc}^{\mathrm{deact}}) \neq \mathbf 0\}$
        \State $M_{\loc} \gets M_{\ell+1}(\mathcal{D}_{\loc}, \mathcal{C}_{\loc}^{\mathrm{deact}})$
       
        \Comment{Step 3: Counter-corrections \eqref{eq:exact_cc}--\eqref{eq:exact_cd}}
        \State $\hat{c}_C \gets c(\mathcal{C}_{\loc}^{\mathrm{coarse}}) - C_\ell^{\ell+1}(\mathcal{C}_{\loc}^{\mathrm{coarse}}, \, \mathcal{A}_{\loc}) \, \hat{c}_A$
        \If{$\mathcal{D}_{\loc} \neq \varnothing$}
        \State $\hat{c}_D \gets c(\mathcal{D}_{\loc}) - M_{\loc} \, C_\ell^{\ell+1}(\mathcal{C}_{\loc}^{\mathrm{deact}}, \, \mathcal{A}_{\loc}) \, \hat{c}_A$
        \EndIf
 
        \Comment{Step 4: Update of the mesh and the basis}
        \State $(\hat{\mathcal{Q}}, \hat{\mathcal{H}}) \gets \text{\bf coarsen\_mesh\_and\_space}(\mathcal{Q}, \mathcal{H}, \ell, \mathcal{Q}_{\loc})$
 
        \Comment{Step 5: Assembly of the coarse coefficient vector}
        \State Initialize a coarse coefficient vector $\hat{c}$ associated with $\hat{\mathcal{H}}$.
        \State $\mathcal{H}_{\mathrm{rem}} \gets\mathcal{H} \setminus (\mathcal{C}_{\loc}^{\fineonly} \cup \mathcal{C}_{\loc}^{\mathrm{coarse}} \cup \mathcal{D}_{\loc})$ \Comment{Same coefficient in both bases} 
        \Statex \textit{Note: by \eqref{eq:hierarchical_basis_disjoint_partitions}, the collections $\mathcal{A}_{\loc}$, $\mathcal{C}_{\loc}^{\mathrm{coarse}}$, $\mathcal{D}_{\loc}$ and $\mathcal{H}_{\mathrm{rem}}$ form a partition of the coarse basis $\hat{\mathcal{H}}$.}
 
        \State $\hat{c}(\mathcal{A}_{\loc}) \gets \hat{c}_A$ \Comment{Insert ancestor coefficients}
        \State $\hat{c}(\mathcal{C}_{\loc}^{\mathrm{coarse}}) \gets \hat{c}_C$ \Comment{Insert corrected coefficients of $\mathcal{C}_{\loc}^{\mathrm{coarse}}$}
        \If{$\mathcal{D}_{\loc} \neq \varnothing$}
        \State $\hat{c}(\mathcal{D}_{\loc}) \gets \hat{c}_D$ \Comment{Insert corrected coefficients of $\mathcal{D}_{\loc}$}
        \EndIf
        \State $\hat{c}(\mathcal{H}_{\mathrm{rem}}) \gets c(\mathcal{H}_{\mathrm{rem}})$ \Comment{Copy unchanged coefficients}
 
        \State \textbf{return} $\hat{\mathcal{Q}}, \hat{\mathcal{H}}, \hat{c}$
    \end{algorithmic}
\end{algorithm}
 
The algorithms of this section take the marked cells as given and return the coarsened mesh, basis and coefficients, together with the bound $e_\infty$ for the $L^\infty$-error, obtained by adding the bounds $e_\ell$ of the single steps. The next section addresses the remaining question: how to choose the marked cells so that the total error does not exceed a prescribed tolerance.

\section{Adaptive coarsening with guaranteed error control}
\label{sec:adaptive_coarsening}

In this section we design a procedure that, for a user-prescribed tolerance $\tau>0$, returns a coarsened mesh together with a coefficient vector whose associated spline is guaranteed to stay within $\tau$ of the original one in the $L^\infty$-norm.

For an arbitrary configuration of marked cells, the error bound of Proposition~\ref{prop:linf_bound} cannot be predicted before it is computed, cf.\ Remark~\ref{rem:l2_vs_linf}, so the tolerance can only be enforced \emph{a posteriori}. This is what the separation between Algorithm~\ref{alg:estimate_error} and Algorithm~\ref{alg:coarsen_at_level} makes possible: the bound is available \emph{before} any topological modification, and an unsatisfactory candidate set can therefore be reduced at no cost other than a new thin least-squares solve. We refer to this correction of the marking as an algebraic \emph{rollback}.

Since a complete sibling group is uniquely identified by its parent cell, we describe the candidates through a set $\hat\QQ_{\loc}\subset\GG_\ell$ of parent cells, as in Section~\ref{S:local_coarsening_problem}, the marked cells being always $\QQ_{\loc}=\children{\hat\QQ_{\loc}}$.

\subsection{Coarsenable cells, error budget and rollback}
\label{S:budget_and_rollback}

\paragraph{Coarsenable cells.} For each level $\ell\in\{0,\dots,L-2\}$ we define the set of \emph{coarsenable cells} of level $\ell$ as
\begin{equation}\label{eq:coarsenable_cells}
    \GG_\ell^{\co}:=\big\{Q\in\GG_\ell\ :\ \children{Q}\subset\QQ\cap\GG_{\ell+1}\big\},
\end{equation}
that is, the cells of level $\ell$ all of whose children are active cells of the hierarchical mesh. Any subset $\hat\QQ_{\loc}\subset\GG_\ell^{\co}$ is an admissible candidate set: the region
\begin{equation}\label{eq:candidate_region}
    \Omega_{\loc}=\bigcup_{Q\in\hat\QQ_{\loc}}Q
\end{equation}
and the marked cells $\QQ_{\loc}=\children{\hat\QQ_{\loc}}$ satisfy~\eqref{eq:sibling_groups_coarsening} by construction, so that no geometric filtering of the kind described in Section~\ref{S:local_coarsening_problem} is required inside the loop.

Initializing the candidate set with the whole set $\GG_\ell^{\co}$, as we do below, amounts to requesting the largest compression that the prescribed tolerance allows. If the coarsening is instead driven by an external marking---for instance the complement of the region marked for refinement in an adaptive time-stepping loop---it suffices to intersect $\GG_\ell^{\co}$ with the corresponding set of cells; nothing else in the algorithm changes. In either case the step remains local: only B-splines interacting with the candidate region enter the computation.

\paragraph{Error budget.} The global tolerance is administered through a scalar budget $\tau_{\mathrm{bud}}$, initialized as $\tau$ and decreased by the exact amount committed at each level. At level $\ell$, the coarsening is accepted only if
\begin{equation}\label{eq:level_admissibility}
    e\ \le\ \tau_{\mathrm{lev}}:=\mu_\ell\,\tau_{\mathrm{bud}},
    \qquad \mu_\ell\in(0,1],
\end{equation}
where $e$ is the bound returned by Algorithm~\ref{alg:estimate_error} for the current candidate set; once the test is passed, that value becomes the error $e_\ell$ committed at this level. The parameters $\mu_\ell$ let the user decide how much of the remaining budget may be consumed at a single level: the choice $\mu_\ell\equiv1$ spends the budget greedily, from the finest transition downwards, whereas $\mu_\ell<1$ reserves part of it for the coarser levels, where a successful coarsening step eliminates larger cells; see Remark~\ref{rem:budget_distribution}.

\paragraph{Rollback.} Assume that the candidate set $\hat\QQ_{\loc}$ has been rejected, i.e.\ that the residual $r_{\loc}$ returned by Algorithm~\ref{alg:estimate_error} satisfies $e=\|r_{\loc}\|_{\infty}>\tau_{\mathrm{lev}}$. We define the set of \emph{violating functions} as
\begin{equation}\label{eq:violating_functions}
    \mathcal C_{\loc}^{\mathrm{viol}}:=
    \begin{cases}
        \big\{\beta\in\mathcal C_{\loc}^{\fineonly}:\ |(r_{\loc})_\beta|=e\big\},
         & \text{(sequential rollback),} \\[4pt]
        \big\{\beta\in\mathcal C_{\loc}^{\fineonly}:\ |(r_{\loc})_\beta|>\tau_{\mathrm{lev}}\big\},
         & \text{(simultaneous rollback),}
    \end{cases}
\end{equation}
and the corresponding set of \emph{violating cells} as
\begin{equation}\label{eq:violating_cells}
    \hat\QQ_{\loc}^{\mathrm{viol}}:=\big\{Q\in\hat\QQ_{\loc}:\
    |\supp(\beta)\cap Q|>0\ \text{ for some }\beta\in\mathcal C_{\loc}^{\mathrm{viol}}\big\}.
\end{equation}
The candidate set is then updated as $\hat\QQ_{\loc}\gets\hat\QQ_{\loc}\setminus\hat\QQ_{\loc}^{\mathrm{viol}}$ and the local problem is solved again. The sequential variant removes only the cells interacting with the largest residual component, and thus discards as few candidates as possible per iteration; the simultaneous variant discards every cell interacting with a violating component at once, and thus needs fewer least-squares solves.

The violating cells~\eqref{eq:violating_cells} comprise \emph{all} the candidate cells met by a violating function, and not only some of them. This choice is dictated by definition~\eqref{eq:def_fineonly}: a function belongs to $\mathcal C_{\loc}^{\fineonly}$ as long as its support overlaps $\Omega_{\loc}$. Discarding all those cells therefore guarantees that a violating function is not removed at the next iteration: its coefficient is kept, and its component leaves the residual. If a single cell $Q$ with $|\supp(\beta)\cap Q|>0$ were kept, $\beta$ would still be removed, and nothing would guarantee that its residual component decreases.

\subsection{The adaptive algorithm}

The structure of Algorithm~\ref{alg:adaptive_coarsen} mirrors that of Algorithm~\ref{alg:coarsen_multilevel}: the levels are swept from the finest transition $\ell=L-2$ down to $\ell=0$, and each level is processed by the two routines of Section~\ref{sec:algorithmic_implementation}. The difference is that the marking is no longer given as input, but determined by the tolerance through the inner loop of Step~2. The set $\GG_\ell^{\co}$ is computed on the current mesh, so that the coarsening of the finer levels may provide new candidates at level~$\ell$. Note also that exiting the inner loop with $\hat\QQ_{\loc}=\varnothing$ is admissible: the level is then left unchanged, and neither the mesh nor the budget is modified.

\begin{algorithm}[h!tb]
    \setstretch{1.25}
    \caption{\bf adaptive\_coarsen\_multilevel}\label{alg:adaptive_coarsen}
    \begin{algorithmic}[1]
        \Statex \textbf{Input:} Hierarchical mesh $\mathcal{Q}$, hierarchical basis $\mathcal{H}$, coefficient vector~$c$, depth $L$ of the hierarchy, tolerance $\tau$, budget fractions $\mu_\ell\in(0,1]$, and a rollback strategy as in \eqref{eq:violating_functions}.
        \Statex \textbf{Output:} Coarsened mesh $\mathcal{Q}$, coarsened basis $\mathcal{H}$, updated coefficient vector $c$, and committed error bound $e_\infty\le\tau$.

        \State $\tau_{\mathrm{bud}} \gets \tau$, \quad $e_\infty\gets0$ \Comment{Initialize the error budget}

        \For{$\ell = L-2$ \textbf{down to} $0$}

        \Comment{Step 1: every coarsenable cell of level $\ell$ is a candidate}
        \State $\hat{\mathcal{Q}}_{\loc} \gets \GG_\ell^{\co}$ as in \eqref{eq:coarsenable_cells}, \quad $\tau_{\mathrm{lev}}\gets\mu_\ell\,\tau_{\mathrm{bud}}$
        \State $\mathrm{accepted} \gets \textbf{false}$, \quad $e_\ell \gets 0$ \Comment{No error is committed unless Step 3 runs}

        \Comment{Step 2: tentative coarsening with algebraic rollback of candidate cells}
        \While{\textbf{not} $\mathrm{accepted}$ \textbf{and} $\hat{\mathcal{Q}}_{\loc} \neq \varnothing$}
        \State $(e, \mathcal{L}_{\loc}, \hat{c}_A, r_{\loc}) \gets \text{\bf local\_coarsening\_residual}(\mathcal{Q}, \mathcal{H}, c, \children{\hat{\mathcal{Q}}_{\loc}}, \ell)$
        \If{$e \leq \tau_{\mathrm{lev}}$}
        \State $\mathrm{accepted} \gets \textbf{true}$, \quad $e_\ell \gets e$
        \Else
        \State Compute $\mathcal C_{\loc}^{\mathrm{viol}}$ by \eqref{eq:violating_functions} and $\hat{\mathcal{Q}}_{\loc}^{\mathrm{viol}}$ by \eqref{eq:violating_cells}
        \State $\hat{\mathcal{Q}}_{\loc} \gets \hat{\mathcal{Q}}_{\loc} \setminus \hat{\mathcal{Q}}_{\loc}^{\mathrm{viol}}$ \Comment{Roll back the violating cells}
        \EndIf
        \EndWhile

        \Comment{Step 3: commit the accepted configuration using the data already computed}
        \If{$\mathrm{accepted}$}
        \State $(\mathcal{Q}, \mathcal{H}, c) \gets \text{\bf coarsen\_at\_level}(\mathcal{Q}, \mathcal{H}, c, \children{\hat{\mathcal{Q}}_{\loc}}, \ell, \mathcal{L}_{\loc}, \hat{c}_A)$
        \State $\tau_{\mathrm{bud}} \gets \tau_{\mathrm{bud}} - e_{\ell}$, \quad $e_\infty\gets e_\infty+e_\ell$ \Comment{Discount the committed error}
        \EndIf
        \EndFor
        \State \textbf{return} $\mathcal{Q}$, $\mathcal{H}$, $c$, $e_\infty$
    \end{algorithmic}
\end{algorithm}

\begin{lemma}[Termination]\label{lem:termination}
    For each level $\ell$, the inner loop of Algorithm~\ref{alg:adaptive_coarsen} performs at most $\#\GG_\ell^{\co}$ iterations, where $\GG_\ell^{\co}$ is computed when level $\ell$ is processed. In particular, the algorithm terminates after at most $\sum_{\ell=0}^{L-2}\#\GG_\ell^{\co}$ calls to Algorithm~\ref{alg:estimate_error}.
\end{lemma}

\begin{proof}
    We first note that the budget never becomes negative: it starts from $\tau>0$, and by~\eqref{eq:level_admissibility} every committed error satisfies $e_\ell\le\mu_\ell\,\tau_{\mathrm{bud}}\le\tau_{\mathrm{bud}}$. Hence $\tau_{\mathrm{lev}}=\mu_\ell\,\tau_{\mathrm{bud}}\ge0$ at every level.
     It suffices to show that every rejected iteration discards at least one cell from $\hat\QQ_{\loc}$. Assume $e>\tau_{\mathrm{lev}}$. Then $e>0$, so that $r_{\loc}$ has at least one component, and since $e=\|r_{\loc}\|_{\infty}$ is attained at some component, in both cases of~\eqref{eq:violating_functions} there exists $\beta\in\mathcal C_{\loc}^{\mathrm{viol}}$. By~\eqref{eq:def_fineonly}, such a $\beta$ satisfies $|\supp(\beta)\cap\Omega_{\loc}|>0$, and $\Omega_{\loc}$ is the union of the finitely many candidate cells~\eqref{eq:candidate_region}; hence $|\supp(\beta)\cap Q|>0$ for at least one $Q\in\hat\QQ_{\loc}$, and $\hat\QQ_{\loc}^{\mathrm{viol}}\neq\varnothing$. Consequently, each rejected iteration decreases $\#\hat\QQ_{\loc}$ by at least one. Since the loop starts with $\#\GG_\ell^{\co}$ candidate cells and every iteration is performed on a nonempty candidate set, at most $\#\GG_\ell^{\co}-1$ rejected iterations precede the last one, and the loop performs at most $\#\GG_\ell^{\co}$ iterations.
\end{proof}

Let us stress that the rollback is a heuristic. Discarding cells modifies both $\mathcal A_{\loc}$ and $\mathcal C_{\loc}^{\fineonly}$, and thus replaces the local problem by a different one, so that the bound $e$ need not decrease monotonically along the loop, and the accepted candidate set is not claimed to be of maximal cardinality among those satisfying~\eqref{eq:level_admissibility}; nor is the resulting compression claimed to grow monotonically with~$\tau$. None of this affects the guarantees of the algorithm: the loop stops, by Lemma~\ref{lem:termination}, and whatever is committed respects the prescribed tolerance, as the next theorem shows.

\begin{theorem}[Success of the adaptive coarsening algorithm]\label{thm:adaptive_success}
    Let $s\in\Span\HH$ be the spline with coefficient vector $c$, and let $\hat s\in\Span\hat\HH$ be the spline returned by Algorithm~\ref{alg:adaptive_coarsen} for an input tolerance $\tau>0$, any budget fractions $\mu_\ell\in(0,1]$ and either rollback strategy, where $\hat\HH$ denotes the basis it returns and $e_\infty$ the committed bound. Then
    \begin{equation*}
        \|s-\hat s\|_{L^\infty(\Omega)}\ \le\ e_\infty\ \le\ \tau .
    \end{equation*}
\end{theorem}

\begin{proof}
    Let $s^{(L-1)}:=s$ and, for $\ell=L-2,\dots,0$, let $s^{(\ell)}$ denote the spline after level $\ell$ has been processed, so that $\hat s=s^{(0)}$. If no coarsening is committed at level $\ell$, then $s^{(\ell)}=s^{(\ell+1)}$ and $e_\ell=0$. Otherwise, $s^{(\ell)}$ is obtained from $s^{(\ell+1)}$ by Algorithm~\ref{alg:coarsen_at_level}, with the local sets and the coefficients $\hat c_A$ that Algorithm~\ref{alg:estimate_error} computed for the same spline and the same candidate set, so that Proposition~\ref{prop:linf_bound} gives $\|s^{(\ell+1)}-s^{(\ell)}\|_{L^\infty(\Omega)}\le e_\ell$. In both cases, the triangle inequality yields
    \begin{equation}\label{eq:telescoping_bound}
        \|s-\hat s\|_{L^\infty(\Omega)}
        \ \le\ \sum_{\ell=0}^{L-2}\|s^{(\ell+1)}-s^{(\ell)}\|_{L^\infty(\Omega)}
        \ \le\ \sum_{\ell=0}^{L-2}e_\ell\ =\ e_\infty .
    \end{equation}
    It remains to check that $e_\infty\le\tau$. Let $\tau^{(L-1)}:=\tau$ and, for $\ell=L-2,\dots,0$, let $\tau^{(\ell)}$ denote the value of $\tau_{\mathrm{bud}}$ after level $\ell$ has been processed, so that
    \begin{equation}\label{eq:budget_recursion}
        \tau^{(\ell)}=\tau^{(\ell+1)}-e_\ell,\qquad \ell=L-2,\dots,0 .
    \end{equation}
    Telescoping~\eqref{eq:budget_recursion} we obtain
    \begin{equation*}
        e_\infty=\sum_{\ell=0}^{L-2}e_\ell
        =\sum_{\ell=0}^{L-2}\big(\tau^{(\ell+1)}-\tau^{(\ell)}\big)
        =\tau-\tau^{(0)},
    \end{equation*}
    so that it suffices to show that $\tau^{(0)}\ge0$. By the acceptance test~\eqref{eq:level_admissibility}, $0\le e_\ell\le\mu_\ell\,\tau^{(\ell+1)}$, an inequality that holds trivially, with $e_\ell=0$, if no coarsening is committed at level $\ell$. Hence $\tau^{(\ell)}\ge(1-\mu_\ell)\,\tau^{(\ell+1)}$, and an induction starting from $\tau^{(L-1)}=\tau>0$ gives $\tau^{(\ell)}\ge0$ for every $\ell$, in particular for $\ell=0$. Together with~\eqref{eq:telescoping_bound}, this completes the proof.
\end{proof}

We remark that the bound of Theorem~\ref{thm:adaptive_success} involves no unknown constant.  Moreover, the quantity $e_\infty$ returned by the algorithm is itself a computable bound, smaller than $\tau$ whenever the budget is not exhausted. In a space-time adaptive loop, this makes it possible to choose $\tau$ as a prescribed fraction of the accuracy demanded from the discretization at the current time step, cf.~\cite{GGMV2019}.

\begin{remark}[Distribution of the budget across the levels]\label{rem:budget_distribution}
    Any sequence $\mu_\ell\in(0,1]$ is admissible in Theorem~\ref{thm:adaptive_success}. The greedy choice $\mu_\ell\equiv1$ allows the finest transition, which is processed first, to consume the whole tolerance, and may therefore leave nothing for the coarser ones, even when these could remove refinement while committing a smaller error. The equidistributed choice $\mu_\ell=1/(\ell+1)$ splits the remaining budget among the $\ell+1$ transitions still to be processed; in particular $\mu_0=1$, so that the coarsest transition may spend whatever has been left unused. The two choices are compared in Section~\ref{S:num_moving_front}, where the equidistributed one removes at least as much refinement as the greedy one while committing a smaller error.
\end{remark}

\begin{remark}[Cost of the adaptive loop]\label{rem:adaptive_cost}
    By Lemma~\ref{lem:termination}, the overhead of the error control amounts to at most $\#\GG_\ell^{\co}$ evaluations of Algorithm~\ref{alg:estimate_error} per level, each of which classifies the local functions and solves a thin overdetermined system as in Remark~\ref{rem:complexity}. The size of these systems never increases along the loop, since discarding cells can only remove functions from $\mathcal A_{\loc}$ and from $\mathcal C_{\loc}^{\fineonly}$. No topological update, and hence no reconstruction of the hierarchical data structures, is ever performed on a rejected configuration: Algorithm~\ref{alg:coarsen_at_level} is called at most once per level, with the classification $\mathcal L_{\loc}$ and the coefficients $\hat c_A$ already computed during the last, accepted, evaluation.
\end{remark}

\begin{remark}[Degenerate configurations]\label{rem:empty_ancestors}
   Some removed functions may not be children of any reactivated one, cf.\ Remark~\ref{rem:no_ancestors}. For such a function $\beta$ the corresponding row of $K_F$ vanishes, so that nothing of it can be recovered and its residual component is its coefficient~$c_\beta$. The adaptive algorithm then behaves as it should: since a configuration is accepted only if every component of $r_{\loc}$ is at most $\tau_{\mathrm{lev}}$, such a function is removed only if $|c_\beta|$ is already that small, and it is then charged in full to the budget; otherwise the cells it overlaps are eventually rolled back.
\end{remark}

\begin{remark}[A localized error budget]\label{rem:local_budget}
    The scalar budget $\tau_{\mathrm{bud}}$ adds up the bounds of all the levels, even when the corresponding perturbations lie in disjoint parts of $\Omega$, and the bound becomes markedly less sharp when several levels commit comparable errors, see Section~\ref{S:num_moving_front}; part of the tolerance is then consumed by errors that do not add up in the actual error. By Proposition~\ref{prop:linf_bound}, however, the coarsening step at level $\ell$ modifies the spline only on the union of the supports of the removed functions. This suggests keeping a separate budget on each cell of a partition of $\Omega$ that does not change during the algorithm, such as the initial mesh $\GG_0$; we say that a step affects a cell of this partition if the cell overlaps that union. More precisely, each cell $Q$ of the partition receives its own budget $\tau_{\mathrm{bud}}(Q)$, initialized to $\tau$. At level $\ell$, the step is accepted if its bound is at most $\mu_\ell$ times the smallest budget among the cells it affects, that is, the most restrictive budget in the region where it modifies the spline; if it is accepted, $e_\ell$ is discounted from the budget of each of those cells, and the other budgets are left unchanged. The guarantee is preserved: on a cell $Q$, the spline changes only in the steps that affect $Q$, and each of them has discounted its bound from $\tau_{\mathrm{bud}}(Q)$; on $Q$ the error is therefore bounded by the sum of those bounds, and the argument of Theorem~\ref{thm:adaptive_success}, applied to $\tau_{\mathrm{bud}}(Q)$ instead of the scalar budget, shows that this sum does not exceed $\tau$. The gain is that steps acting on separate regions no longer compete for the same tolerance: the scalar budget discounts every step, whereas $\tau_{\mathrm{bud}}(Q)$ discounts only those affecting $Q$, so that, for a given sequence of committed steps, $\tau_{\mathrm{bud}}(Q)$ is never smaller than the scalar budget and the acceptance test is never more restrictive; more refinement may therefore be removed for the same~$\tau$. We leave the numerical study of this variant for future work.
\end{remark}

\section{Coarsening for truncated hierarchical B-splines}
\label{sec:thb}

The truncated hierarchical B-spline basis of~\cite{GJS12} spans the same space as $\HH$ while being a partition of unity, and is the basis of choice in most implementations. We show in this section that the local framework of Section~\ref{sec:coarsening} carries over to it, and that truncation in fact simplifies the coarsening step.

Truncation removes from each coarse function the components that lie inside the refined region, and a coarsening step gives some of them back. Thus, a reactivated ancestor is expressed \emph{exclusively} in terms of the functions that the step removes: what truncation had discarded from it is precisely what, in the HB setting, had to be counter-corrected on $\mathcal C_{\loc}^{\mathrm{coarse}}$ and $\mathcal D_{\loc}$. In exchange, the coarse functions of levels $\le\ell$ that survive the step are \emph{different functions} in the two bases, since part of what was truncated is returned to them. We show below that the local problem reduces to a least-squares problem with the same matrix $K_F$ as in Section~\ref{sec:coarsening}, but with a modified right-hand side.

\paragraph{The truncated basis.} Let $\boldsymbol\Omega_L$ be a hierarchy of subdomains as in Section~\ref{sec:hierarchical_bsplines}. For $1\le m\le L-1$ and $f\in\mathcal S_{m-1}$, the \emph{truncation} of $f$ at level $m$ discards, from the representation of $f$ in $\mathcal B_m$, the contributions of the B-splines whose support lies in $\Omega_m$,
\begin{equation}\label{eq:trunc_def}
	\operatorname{trunc}_m(f):=\sum_{\substack{\gamma\in\mathcal B_m \\ \supp\gamma\not\subset\Omega_m}} c_\gamma(f)\,\gamma,
	\qquad\text{where}\qquad
	f=\sum_{\gamma\in\mathcal B_m}c_\gamma(f)\,\gamma .
\end{equation}
The discarded functions are those entering the hierarchical basis at step $m$, that is, $\{\gamma\in\mathcal B_m:\supp\gamma\subset\Omega_m\}=\mathcal F_m^A\,\dot\cup\,\mathcal F_m^D$; setting to zero the corresponding rows of $C_{m-1}^{m}$ yields the truncated subdivision matrix $C_{m-1}^{m,\mathrm{tr}}$, which maps the coefficients of $f$ in $\mathcal B_{m-1}$ to those of $\operatorname{trunc}_m(f)$ in $\mathcal B_m$. For $\beta\in\mathcal B_k$, $0\le k\le L-1$, applying~\eqref{eq:trunc_def} successively at the levels above~$k$ defines the \emph{partial truncations} of~$\beta$,
\begin{equation}\label{eq:partial_truncation}
	T^{(k)}\beta:=\beta,
	\qquad
	T^{(m)}\beta:=\operatorname{trunc}_m\big(T^{(m-1)}\beta\big)\ \ (k<m\le L-1),
	\qquad
	T\beta:=T^{(L-1)}\beta.
\end{equation}
Each step is well defined because $T^{(m-1)}\beta\in\mathcal S_{m-1}$, and $T\beta=\beta$ when $k=L-1$. Applied to the active B-splines, this construction produces the THB-spline basis $\Tau:=\{T\beta:\beta\in\mathcal F^A\}$ of~\cite{GJS12}, which satisfies $\Span\Tau=\Span\HH$. Notice that $T\beta$ depends not only on $\beta$ but also on the subdomains $\Omega_{k+1},\dots,\Omega_{L-1}$: a change of the mesh may alter $T\beta$ even when $\beta$ remains active. This property has no counterpart in the HB setting, and it is the reason why the analysis of Section~\ref{sec:coarsening} does not apply directly to THB-splines.
 
In Section~\ref{sec:hierarchical_bsplines}, $\mathcal F^A=\HH$ and an active B-spline is at once a label and an element of the basis. Here the two roles separate, since the elements of $\Tau$ are labelled by the active B-splines without being, in general, those B-splines. We therefore read every index set of the paper as a set of B-splines, so that the coefficient vectors and the submatrix notation $M(\mathcal X,\mathcal Y)$, $v(\mathcal X)$ of Section~\ref{sec:hierarchical_bsplines} apply unchanged. Besides changing the active set, a coarsening step may thus change the function attached to an index that it keeps: $T\beta$ in the fine configuration, $\hat T\beta$ in the coarse one, where $\hat T$ is built as in~\eqref{eq:partial_truncation} from the coarse hierarchy~\eqref{eq:coarse_hierarchy}.

\paragraph{The subdivision operators.} Following~\cite{GV18}, the partial truncations of all the active functions up to level $\ell$ are collected in a single matrix: writing $\mathcal F^A_{0:\ell}:=\bigcup_{k\le\ell}\mathcal F^A_k$, the matrix $C^{\mathrm{sub}}_\ell\in\RR^{N_\ell\times\#\mathcal F^A_{0:\ell}}$, $0\le\ell\le L-1$, is defined by
\begin{equation}\label{eq:csub_def}
	T^{(\ell)}\beta=\sum_{\gamma\in\mathcal B_\ell}C^{\mathrm{sub}}_\ell(\gamma,\beta)\,\gamma,
	\qquad \beta\in\mathcal F^A_{0:\ell}.
\end{equation}
By~\eqref{eq:partial_truncation}, $T^{(\ell)}\beta=\beta$ for $\beta\in\mathcal F^A_\ell$, while for $\beta\in\mathcal F^A_{0:\ell}$ the function $T^{(\ell+1)}\beta=\operatorname{trunc}_{\ell+1}\big(T^{(\ell)}\beta\big)$ has coefficients $C_\ell^{\ell+1,\mathrm{tr}}\,C^{\mathrm{sub}}_\ell(\cdot,\beta)$ in $\mathcal B_{\ell+1}$. Ordering the columns by level, the matrices are therefore given by the recursion
\begin{equation}\label{eq:csub_recursion}
	C^{\mathrm{sub}}_0=J_0,
	\qquad
	C^{\mathrm{sub}}_{\ell+1}=\big[\,C_\ell^{\ell+1,\mathrm{tr}}\,C^{\mathrm{sub}}_\ell,\ J_{\ell+1}\,\big],
	\qquad 0\le\ell\le L-2,
\end{equation}
where $J_\ell:=I(\mathcal B_\ell,\mathcal F^A_\ell)\in\RR^{N_\ell\times N_\ell^A}$ consists of the columns of the identity matrix of order $N_\ell$ indexed by the active B-splines of level~$\ell$.
 The entries of $C^{\mathrm{sub}}_\ell$ are non-negative, since the recursion only multiplies non-negative subdivision weights and sets entries to zero. For $\ell=L-1$,~\eqref{eq:csub_def} expresses every function of $\Tau$ in $\mathcal B_{L-1}$, so that the THB-splines are non-negative. In GeoPDEs~\cite{GV18} these matrices are assembled by the function \texttt{hspace\_subdivision\_matrix} and stored in the property \texttt{Csub} of the \texttt{hierarchical\_space} class. Since the recursion only involves subdivision matrices, $T^{(m)}\beta$ is a linear combination of descendants of $\beta$; in particular, $\supp T\beta\subseteq\supp\beta$.
 
\subsection{The refinement matrix between truncated bases}
\label{S:thb_block_structure}
 
The setting is that of Section~\ref{S:local_coarsening_problem}: a set of active cells $\QQ_{\loc}=\children{\hat\QQ_{\loc}}\subset \GG_{\ell+1}\cap\QQ$ is marked for coarsening, $\Omega_{\loc}$ is the region~\eqref{eq:local_region} and $\hat{\boldsymbol\Omega}_L$ is the coarse hierarchy~\eqref{eq:coarse_hierarchy}, whose truncated basis we denote by $\hat\Tau$. Since the status of a B-spline depends only on the hierarchy and not on the truncation, the sets $\mathcal A_{\loc}$ and $\mathcal C_{\loc}^{\fineonly}$ of~\eqref{eq:def_Aloc}--\eqref{eq:def_fineonly} are those of Section~\ref{S:local_block_structure}, the active sets are related by
\begin{equation}\label{eq:thb_active_update}
	\hat{\mathcal F}^A=\big(\mathcal F^A\setminus\mathcal C_{\loc}^{\fineonly}\big)\ \dot\cup\ \mathcal A_{\loc},
\end{equation}
and the dimension update~\eqref{eq:coarsened_dimension_update} is unchanged.

Since $\Span\hat\Tau=\Span\hat\HH\subseteq\Span\HH=\Span\Tau$, every function of the coarse basis can be uniquely expanded in the fine one. We collect these expansions in the refinement matrix $K^{\mathrm{trunc}}\in\RR^{n\times\hat n}$,
\begin{equation*}
	\hat T\beta=\sum_{\gamma\in\mathcal F^A}K^{\mathrm{trunc}}(\gamma,\beta)\,T\gamma,
	\qquad \beta\in\hat{\mathcal F}^A,
\end{equation*}
whose column $\beta$ holds the coefficients of $\hat T\beta$ in $\Tau$. As in~\eqref{eq:refinement_relation}, it transfers coefficient vectors: a spline with coefficients $\hat c$ in $\hat\Tau$ has coefficients $c=K^{\mathrm{trunc}}\hat c$ in $\Tau$.

By~\eqref{eq:coarse_hierarchy}, the two configurations share every subdomain but $\Omega_{\ell+1}$; their meshes therefore differ only in the marked region, $\hat\QQ=(\QQ\setminus\QQ_{\loc})\cup\hat\QQ_{\loc}$. As for the bases, both are built by~\eqref{eq:partial_truncation}, and the operators $\operatorname{trunc}_m$, each of which depends only on $\Omega_m$, all coincide except $\operatorname{trunc}_{\ell+1}$, which discards the B-splines $\gamma\in\mathcal B_{\ell+1}$ with $\supp\gamma\subset\Omega_{\ell+1}$ in the fine configuration and those with $\supp\gamma\subset\hat\Omega_{\ell+1}$ in the coarse one. As noted after~\eqref{eq:partial_truncation}, the truncation of a level-$k$ function involves only the subdomains $\Omega_{k+1},\dots,\Omega_{L-1}$, which do not include $\Omega_{\ell+1}$ when $k\ge\ell+1$, so that
\begin{equation}\label{eq:thb_fine_levels_unchanged}
    \hat T\beta=T\beta\qquad\text{for every }\beta\in\hat{\mathcal F}^A\cap\mathcal B_k\text{ with }k\ge\ell+1.
\end{equation}
   By~\eqref{eq:thb_active_update}, these indices are also active in the fine configuration, so that $T\beta\in\Tau$ and the column $K^{\mathrm{trunc}}(\cdot,\beta)$ vanishes except for $K^{\mathrm{trunc}}(\beta,\beta)=1$. There remain the columns indexed by the coarse functions of levels $\le\ell$, that is, by $\mathcal F^A_{0:\ell}\cup\mathcal A_{\loc}$, for which the comparison of $\hat T\beta$ with $T\beta$ reduces to comparing the two versions of $\operatorname{trunc}_{\ell+1}$. We do this in two steps, identifying first the B-splines on which they differ.

\begin{lemma}[Effect of coarsening on $\operatorname{trunc}_{\ell+1}$]\label{lem:thb_fineonly_characterization}
    The B-splines discarded by $\operatorname{trunc}_{\ell+1}$ in the fine configuration but not in the coarse one are exactly the fine-only functions, that is,
    \begin{equation}\label{eq:thb_fineonly_characterization}
        \big\{\gamma\in\mathcal B_{\ell+1}:\ \supp\gamma\subset\Omega_{\ell+1}
        \ \text{ and }\ \supp\gamma\not\subset\hat\Omega_{\ell+1}\big\}
        =\mathcal C_{\loc}^{\fineonly} .
    \end{equation}
\end{lemma}
 
\begin{proof}
    Since $\Omega_{\loc}$ is a union of cells of $\GG_{\ell+1}$ and $\hat\Omega_{\ell+1}$ is the closure of $\Omega_{\ell+1}\setminus\Omega_{\loc}$, every cell $Q\in\GG_{\ell+1}$ with $Q\subset\Omega_{\ell+1}$ lies either in $\Omega_{\loc}$ or in $\hat\Omega_{\ell+1}$, but not in both. If $\gamma\in\mathcal B_{\ell+1}$ satisfies $\supp\gamma\subset\Omega_{\ell+1}$, its support is a union of such cells, and therefore $\supp\gamma\not\subset\hat\Omega_{\ell+1}$ if and only if $|\supp\gamma\cap\Omega_{\loc}|>0$.

    Let now $\gamma$ belong to the left-hand side of~\eqref{eq:thb_fineonly_characterization}. Since $\supp\gamma\subset\Omega_{\ell+1}$, $\gamma$ is either active or deactivated at level $\ell+1$. It cannot be deactivated: otherwise $\supp\gamma\subset\Omega_{\ell+2}$, and $\Omega_{\ell+2}$ is a union of cells of $\GG_{\ell+1}$ that are not active, hence not marked, so that $\Omega_{\ell+2}\subset\hat\Omega_{\ell+1}$, contradicting $\supp\gamma\not\subset\hat\Omega_{\ell+1}$. Thus $\gamma\in\mathcal F^A_{\ell+1}$ and $|\supp\gamma\cap\Omega_{\loc}|>0$, that is, $\gamma\in\mathcal C_{\loc}^{\fineonly}$ by~\eqref{eq:def_fineonly}. Conversely, every $\gamma\in\mathcal C_{\loc}^{\fineonly}$ is active of level $\ell+1$, so that $\supp\gamma\subset\Omega_{\ell+1}$, and $|\supp\gamma\cap\Omega_{\loc}|>0$ gives $\supp\gamma\not\subset\hat\Omega_{\ell+1}$.
\end{proof}

Coarsening modifies $\operatorname{trunc}_{\ell+1}$, which could a priori alter the components of the coarse functions of levels $\le\ell$ along any B-spline of level $\ell+1$. Lemma~\ref{lem:thb_fineonly_characterization} shows that the change is local: the two versions differ only on the fine-only functions, whose supports overlap $\Omega_{\loc}$ and which the step removes. It remains to propagate that difference through the truncations of the levels above $\ell+1$, which are common to both configurations.

To this end, we collect in a matrix $E_{\ell+1}$ the coefficients in $\mathcal B_{\ell+1}$ of the partial truncations of level~$\ell$ of the coarse functions of levels $\le\ell$,
\begin{equation}\label{eq:def_E}
    T^{(\ell)}\beta=\sum_{\gamma\in\mathcal B_{\ell+1}}E_{\ell+1}(\gamma,\beta)\,\gamma,
    \qquad \beta\in\mathcal F^A_{0:\ell}\cup\mathcal A_{\loc}.
\end{equation}
By~\eqref{eq:csub_def} and the two-scale relation~\eqref{E:two scale relation}, $E_{\ell+1}(\cdot,\beta)=C_\ell^{\ell+1}\,C^{\mathrm{sub}}_\ell(\cdot,\beta)$ for $\beta\in\mathcal F^A_{0:\ell}$, while $E_{\ell+1}(\cdot,\beta)=C_\ell^{\ell+1}(\cdot,\beta)$ for $\beta\in\mathcal A_{\loc}$, because $T^{(\ell)}\beta=\beta$ for these level-$\ell$ B-splines.

\begin{lemma}[Expansion of the coarse functions]\label{lem:thb_expansion}
   Every coarse function of level at most $\ell$ differs from the function attached to the same index in the fine configuration by a linear combination of the functions $T\gamma$, $\gamma\in\mathcal C_{\loc}^{\fineonly}$:
    \begin{equation}\label{eq:thb_expansion}
        \hat T\beta=T\beta+\sum_{\gamma\in\mathcal C_{\loc}^{\fineonly}}E_{\ell+1}(\gamma,\beta)\,T\gamma,
        \qquad \beta\in\mathcal F^A_{0:\ell}\cup\mathcal A_{\loc}.
    \end{equation}
    Moreover, $T\beta=0$ for every $\beta\in\mathcal A_{\loc}$.
\end{lemma}

\begin{proof}
    Let $\beta\in\mathcal F^A_{0:\ell}\cup\mathcal A_{\loc}$. The operators $\operatorname{trunc}_k$ with $k\le\ell$ are the same in both configurations, so $\hat T^{(\ell)}\beta=T^{(\ell)}\beta$, with coefficients $E_{\ell+1}(\cdot,\beta)$ in $\mathcal B_{\ell+1}$ by~\eqref{eq:def_E}. At level $\ell+1$, the coarse version of $\operatorname{trunc}_{\ell+1}$ keeps all the terms kept by the fine one and, in addition, those along the B-splines of the left-hand side of~\eqref{eq:thb_fineonly_characterization}, which are the fine-only functions by Lemma~\ref{lem:thb_fineonly_characterization}. Hence
    \begin{equation*}
        \hat T^{(\ell+1)}\beta-T^{(\ell+1)}\beta
        =\sum_{\gamma\in\mathcal C_{\loc}^{\fineonly}}E_{\ell+1}(\gamma,\beta)\,\gamma .
    \end{equation*}
    By~\eqref{eq:partial_truncation}, the functions $\hat T\beta$, $T\beta$ and $T\gamma$, $\gamma\in\mathcal B_{\ell+1}$, result from $\hat T^{(\ell+1)}\beta$, $T^{(\ell+1)}\beta$ and $\gamma$, respectively, through the operators $\operatorname{trunc}_{\ell+2},\dots,\operatorname{trunc}_{L-1}$ (there are none if $\ell=L-2$). Since these operators are linear and common to both configurations, applying them to the last identity yields~\eqref{eq:thb_expansion}. Finally, if $\beta\in\mathcal A_{\loc}\subset\mathcal F_\ell^D$, then $\supp\gamma\subset\supp\beta\subset\Omega_{\ell+1}$ for every child $\gamma$ of $\beta$, so that the fine version of $\operatorname{trunc}_{\ell+1}$ discards all the level-$(\ell+1)$ components of $\beta$. Thus $T^{(\ell+1)}\beta=0$, and therefore $T\beta=0$.
\end{proof}

By Lemma~\ref{lem:thb_expansion} and~\eqref{eq:thb_fine_levels_unchanged}, an index $\beta$ that is active in both configurations, that is, $\beta\in\mathcal F^A\setminus\mathcal C_{\loc}^{\fineonly}=\hat{\mathcal F}^A\setminus\mathcal A_{\loc}$ by~\eqref{eq:thb_active_update}, satisfies $\hat T\beta\neq T\beta$ if and only if it has level $\le\ell$ and $E_{\ell+1}(\mathcal C_{\loc}^{\fineonly},\beta)\neq\mathbf 0$. We collect these indices in $\mathcal F_{\loc}^{\mathrm{trunc}}$, and the remaining ones, whose function does not change, in $\mathcal F_{\mathrm{rem}}$:
\begin{equation}\label{eq:def_Ftrunc}
    \mathcal F_{\loc}^{\mathrm{trunc}}:=\big\{\beta\in\mathcal F^A_{0:\ell}:\ E_{\ell+1}\big(\mathcal C_{\loc}^{\fineonly},\beta\big)\neq{\bf 0}\big\},
    \qquad
    \mathcal F_{\mathrm{rem}}:=\mathcal F^A\setminus\big(\mathcal C_{\loc}^{\fineonly}\cup\mathcal F_{\loc}^{\mathrm{trunc}}\big).
\end{equation}
Only the active functions of levels $\le\ell$ whose support overlaps $\Omega_{\loc}$ need to be tested in~\eqref{eq:def_Ftrunc}, because $E_{\ell+1}(\gamma,\beta)\neq0$ implies $\supp\gamma\subset\supp\beta$. With $n_T:=\#\mathcal F_{\loc}^{\mathrm{trunc}}$, we set
\begin{equation}\label{eq:def_KT}
\begin{aligned}
    K_F&:=E_{\ell+1}\big(\mathcal C_{\loc}^{\fineonly},\mathcal A_{\loc}\big)
    =C_\ell^{\ell+1}\big(\mathcal C_{\loc}^{\fineonly},\mathcal A_{\loc}\big)\in\RR^{n_F\times n_A},\\
    K_T&:=E_{\ell+1}\big(\mathcal C_{\loc}^{\fineonly},\mathcal F_{\loc}^{\mathrm{trunc}}\big)\in\RR^{n_F\times n_T}.
\end{aligned}
\end{equation}
In practice, $K_T$ cannot be extracted from the stored matrix $C^{\mathrm{sub}}_{\ell+1}$, since by~\eqref{eq:csub_recursion} its rows indexed by $\mathcal C_{\loc}^{\fineonly}\subset\mathcal F^A_{\ell+1}$ have been set to zero by the truncation; it is computed instead as $K_T=C_\ell^{\ell+1}(\mathcal C_{\loc}^{\fineonly},\cdot)\,C^{\mathrm{sub}}_\ell(\cdot,\mathcal F_{\loc}^{\mathrm{trunc}})$, with the untruncated matrix $C_\ell^{\ell+1}$. 

\begin{proposition}[Block structure of $K^{\mathrm{trunc}}$]\label{prop:thb_block_structure}
    The sets $\mathcal C_{\loc}^{\fineonly}$, $\mathcal F_{\loc}^{\mathrm{trunc}}$, $\mathcal F_{\mathrm{rem}}$ form a partition of $\mathcal F^A$, and $\mathcal A_{\loc}$, $\mathcal F_{\loc}^{\mathrm{trunc}}$, $\mathcal F_{\mathrm{rem}}$ form a partition of $\hat{\mathcal F}^A$. If the rows and the columns of $K^{\mathrm{trunc}}$ are ordered accordingly, then
    \begin{equation}\label{eq:thb_block}
        K^{\mathrm{trunc}}=
        \begin{blockarray}{ccc l}
            \mathcal A_{\loc} & \mathcal F_{\loc}^{\mathrm{trunc}} & \mathcal F_{\mathrm{rem}} & \\
            \begin{block}{(ccc)l}
                K_F & K_T & 0 & \;\mathcal C_{\loc}^{\fineonly} \\
                0 & I & 0 & \;\mathcal F_{\loc}^{\mathrm{trunc}} \\
                0 & 0 & I & \;\mathcal F_{\mathrm{rem}} \\
            \end{block}
        \end{blockarray},
    \end{equation}
    where $K_F$ and $K_T$ are the blocks~\eqref{eq:def_KT} and $I$ denotes identity matrices of the appropriate sizes.
\end{proposition}

\begin{proof}
    The sets $\mathcal F_{\loc}^{\mathrm{trunc}}\subset\mathcal F^A_{0:\ell}$ and $\mathcal C_{\loc}^{\fineonly}\subset\mathcal F^A_{\ell+1}$ are disjoint, so the first partition holds by the definition of $\mathcal F_{\mathrm{rem}}$; the second one follows from~\eqref{eq:thb_active_update}, since $\mathcal A_{\loc}\subset\mathcal F^D_\ell$ is disjoint from $\mathcal F^A$.

    Since $\Tau$ is a basis, the column of $K^{\mathrm{trunc}}$ indexed by $\beta\in\hat{\mathcal F}^A$ is the unique coefficient vector of $\hat T\beta$ in $\Tau$. If $\beta$ has level $\ge\ell+1$, then $\beta\in\mathcal F_{\mathrm{rem}}$ and, by~\eqref{eq:thb_fine_levels_unchanged}, $K^{\mathrm{trunc}}(\cdot,\beta)$ vanishes except for $K^{\mathrm{trunc}}(\beta,\beta)=1$. Otherwise $\beta\in\mathcal F^A_{0:\ell}\cup\mathcal A_{\loc}$, and~\eqref{eq:thb_expansion} gives the column: its entries on the rows of $\mathcal C_{\loc}^{\fineonly}$ are $E_{\ell+1}(\mathcal C_{\loc}^{\fineonly},\beta)$, its entry on the row of $\beta$ is $1$ if $\beta\in\mathcal F^A_{0:\ell}$, and all its other entries vanish, since $T\beta=0$ for $\beta\in\mathcal A_{\loc}$. Reading these columns by blocks of rows gives~\eqref{eq:thb_block}. For $\beta\in\mathcal A_{\loc}$, the column has the entries $K_F(\cdot,\beta)$ on the rows of $\mathcal C_{\loc}^{\fineonly}$ and vanishes elsewhere. For $\beta\in\mathcal F_{\loc}^{\mathrm{trunc}}$, it has the entries $K_T(\cdot,\beta)$ on the rows of $\mathcal C_{\loc}^{\fineonly}$ and, outside them, a single nonzero entry, equal to $1$, on the row of $\beta$. For $\beta\in\mathcal F_{\mathrm{rem}}$, its entries on the rows of $\mathcal C_{\loc}^{\fineonly}$ vanish by~\eqref{eq:def_Ftrunc},    so that the column again vanishes except for the entry $1$ on the row of $\beta$.
\end{proof}

\begin{remark}[Comparison with the HB setting]\label{rem:thb_vs_hb}
    Table~\ref{tab:thb_local_sets} and the block structure~\eqref{eq:thb_block} are the counterparts of Table~\ref{tab:local_sets} and~\eqref{eq:global_refinement_detailed_block}. In the HB setting, reactivating an ancestor adds contributions to the surviving functions of $\mathcal C_{\loc}^{\mathrm{coarse}}$ and $\mathcal D_{\loc}$, which produce the blocks $K_C$ and $K_D$ and require the counter-corrections~\eqref{eq:exact_cc}--\eqref{eq:exact_cd}. With truncation, the column of an ancestor reduces to the single block $K_F$; in exchange, the surviving functions of $\mathcal F_{\loc}^{\mathrm{trunc}}$, whose columns are unit vectors in the HB setting, now carry the block $K_T$. Both features have the same origin: $\operatorname{trunc}_{\ell+1}$ discards the components supported in the refined subdomain, which is $\Omega_{\ell+1}$ in the fine configuration and $\hat\Omega_{\ell+1}$ in the coarse one. In the coarse configuration, the children of an ancestor in $\mathcal C_{\loc}^{\mathrm{coarse}}$ and $\mathcal C_{\loc}^{\mathrm{deact}}$ have their supports in $\hat\Omega_{\ell+1}$ and are discarded, which removes $K_C$ and $K_D$; a surviving coarse function, on the other hand, recovers its components along the fine-only functions, which form $K_T$.

\begin{table}[htb]
    \footnotesize
    \centering
    \renewcommand{\arraystretch}{1.25}
    \begin{tabular}{lccccc}
        \toprule
        \textbf{Set} & \textbf{Level} & \textbf{Coarse status} & \textbf{Fine status} & \textbf{Same function?} & \textbf{Coefficient} \\
        \midrule
        $\mathcal A_{\loc}$ & $\ell$ & active & deactivated & --- & recovered by least squares \\
        $\mathcal C_{\loc}^{\fineonly}$ & $\ell+1$ & passive & active & --- & removed \\
        $\mathcal F_{\loc}^{\mathrm{trunc}}$ & $\le\ell$ & active & active & no & copied \\
        $\mathcal F_{\mathrm{rem}}$ & any & active & active & yes & copied \\
        \bottomrule
    \end{tabular}
    \caption{The local index sets in the truncated setting, their status in each configuration, whether the two configurations attach the same function to them, and the treatment of their coefficients in the coarsening step of Section~\ref{S:thb_local_problem}.}
    \label{tab:thb_local_sets}
\end{table}
\end{remark}

\begin{remark}[Partition of unity]\label{rem:thb_partition_unity}
    Since $\Tau$ and $\hat\Tau$ are partitions of unity, the constant function has all its coefficients equal to one in both bases, and~\eqref{eq:thb_block} yields
    \begin{equation}\label{eq:thb_partition_unity_blocks}
        K_F\mathbf 1_{\mathcal A_{\loc}}
        +K_T\mathbf 1_{\mathcal F_{\loc}^{\mathrm{trunc}}}
        =\mathbf 1_{\mathcal C_{\loc}^{\fineonly}}.
    \end{equation}
    This identity provides a simple test for the implementation of both blocks. It also shows that the surviving functions cannot be ignored. If $\mathcal F_{\loc}^{\mathrm{trunc}}\neq\varnothing$, the columns of $K_T$ are nonzero by~\eqref{eq:def_Ftrunc} and have non-negative entries, so that $K_T\mathbf 1_{\mathcal F_{\loc}^{\mathrm{trunc}}}\neq\mathbf 0$ and hence $K_F\mathbf 1_{\mathcal A_{\loc}}\neq\mathbf 1_{\mathcal C_{\loc}^{\fineonly}}$: the ancestors only recover the part of the constant that is not already carried by the functions of $\mathcal F_{\loc}^{\mathrm{trunc}}$.
\end{remark}

\subsection{The local coarsening problem and the error bound}
\label{S:thb_local_problem}
 
 Let $s=\sum_{\beta\in\mathcal F^A}c_\beta\,T\beta$ be the spline to be coarsened, with coefficient vector $c$ in $\Tau$. A coarsening step should alter $s$ only where the mesh changes; as in~\eqref{eq:admissible_class}, we therefore look for the coarse spline in the class
\begin{equation}\label{eq:thb_admissible_class}
    \hat{\mathcal V}:=\big\{\hat v\in\Span\hat\Tau:\ \text{the coefficients of }\hat v\text{ and }s\text{ in }\Tau\text{ coincide on }\mathcal F^A\setminus\mathcal C_{\loc}^{\fineonly}\big\} .
\end{equation}
The indices of $\mathcal F^A\setminus\mathcal C_{\loc}^{\fineonly}=\mathcal F_{\loc}^{\mathrm{trunc}}\cup\mathcal F_{\mathrm{rem}}$ are those that survive the step, and by the two identity blocks of~\eqref{eq:thb_block} the coefficients of $\hat v$ on them are the same in $\Tau$ and in $\hat\Tau$. Hence an element of $\hat{\mathcal V}$ has the coarse coefficients of $s$ on every surviving index, and is determined by its \emph{ancestor coefficients} $a\in\RR^{n_A}$, that is, by its coefficients in $\hat\Tau$ on $\mathcal A_{\loc}$. We abbreviate
\begin{equation}\label{eq:thb_local_vectors}
    c_F:=c\big(\mathcal C_{\loc}^{\fineonly}\big),
    \qquad
    c_T:=c\big(\mathcal F_{\loc}^{\mathrm{trunc}}\big),
    \qquad
    g_F:=c_F-K_T\,c_T .
\end{equation}

\begin{lemma}[Error representation]\label{lem:thb_localization}
    Let $\hat v\in\hat{\mathcal V}$ have ancestor coefficients $a\in\RR^{n_A}$. Then $s-\hat v$ is a linear combination of the functions removed by the step,
    \begin{equation}\label{eq:thb_difference}
        s-\hat v=\sum_{\gamma\in\mathcal C_{\loc}^{\fineonly}}\big(g_F-K_F a\big)_\gamma\ T\gamma .
    \end{equation}
\end{lemma}

\begin{proof}
    The coefficient vector of $\hat v$ in $\hat\Tau$ is $a$ on $\mathcal A_{\loc}$ and coincides with $c$ on $\mathcal F_{\loc}^{\mathrm{trunc}}\cup\mathcal F_{\mathrm{rem}}$. Multiplying it by $K^{\mathrm{trunc}}$, in the form~\eqref{eq:thb_block}, gives its coefficients in $\Tau$: $K_Fa+K_Tc_T$ on $\mathcal C_{\loc}^{\fineonly}$, and $c$ on $\mathcal F_{\loc}^{\mathrm{trunc}}\cup\mathcal F_{\mathrm{rem}}$. Hence the coefficients of $s-\hat v$ in $\Tau$ vanish outside $\mathcal C_{\loc}^{\fineonly}$, and on $\mathcal C_{\loc}^{\fineonly}$ they are $c_F-K_Tc_T-K_Fa=g_F-K_Fa$.
\end{proof}

By Lemma~\ref{lem:thb_localization}, the whole step is governed by the system $K_Fa\approx g_F$. The ancestors must reproduce not the removed coefficients $c_F$ themselves, but what remains of them after subtracting $K_Tc_T$, the part already recovered by the surviving coarse functions, as observed for the constant function in Remark~\ref{rem:thb_partition_unity}. The matrix $K_F=C_\ell^{\ell+1}(\mathcal C_{\loc}^{\fineonly},\mathcal A_{\loc})$ is the same as in the HB setting, cf.~\eqref{eq:coarsening_local_blocks}, so that Lemma~\ref{lemma:kf_rank} applies: $K_F$ has full column rank, hence $n_A\le n_F$ and the system has a unique least-squares solution. We define the ancestor coefficients of the coarse spline and the residual by
\begin{equation}\label{eq:thb_ls_problem}
    \hat c_A:=\argmin_{a\in\RR^{n_A}}\big\|g_F-K_Fa\big\|_2=K_F^\dagger g_F,
    \qquad
    r_{\loc}:=g_F-K_F\hat c_A,
\end{equation}
and we let $\Pi^{\mathrm{trunc}}:\Span\Tau\to\Span\hat\Tau$ denote the resulting coarsening operator, which maps $s$ to the element $\hat s$ of $\hat{\mathcal V}$ with ancestor coefficients $\hat c_A$, that is, to the spline whose coefficient vector $\hat c$ in $\hat\Tau$ is
\begin{equation}\label{eq:thb_coarse_coefficients}
    \hat c(\mathcal A_{\loc})=\hat c_A,
    \qquad
    \hat c\big(\mathcal F_{\loc}^{\mathrm{trunc}}\cup\mathcal F_{\mathrm{rem}}\big)=c\big(\mathcal F_{\loc}^{\mathrm{trunc}}\cup\mathcal F_{\mathrm{rem}}\big).
\end{equation}
Since $\hat c$ depends linearly on $c$, $\Pi^{\mathrm{trunc}}$ is linear. The block structure~\eqref{eq:thb_block} yields at once the counterpart of Corollary~\ref{cor:exactness}.

\begin{proposition}[Left inverse and projection]\label{prop:thb_exactness}
    The map $c\mapsto\hat c$ is a left inverse of $K^{\mathrm{trunc}}$. Consequently, if $s\in\Span\hat\Tau$, then $\Pi^{\mathrm{trunc}}s=s$ and the residual $r_{\loc}$ of~\eqref{eq:thb_ls_problem} vanishes; in particular, $\Pi^{\mathrm{trunc}}$ is a projection onto $\Span\hat\Tau$.
\end{proposition}

\begin{proof}
       By~\eqref{eq:thb_local_vectors}, \eqref{eq:thb_ls_problem} and~\eqref{eq:thb_coarse_coefficients}, $\hat c=Y^{\mathrm{trunc}}c$, where, with blocks ordered as in~\eqref{eq:thb_block},
    \[
        Y^{\mathrm{trunc}}:=\begin{pmatrix} K_F^\dagger & -K_F^\dagger K_T & 0 \\ 0 & I & 0 \\ 0 & 0 & I \end{pmatrix}.
    \]
    Since $K_F^\dagger K_F=I$, block multiplication with~\eqref{eq:thb_block} gives $Y^{\mathrm{trunc}}K^{\mathrm{trunc}}=I$. Now let $s\in\Span\hat\Tau$, with coefficient vector $\hat v$ in $\hat\Tau$. Then $c=K^{\mathrm{trunc}}\hat v$, so that $\hat c=Y^{\mathrm{trunc}}K^{\mathrm{trunc}}\hat v=\hat v$, that is, $\Pi^{\mathrm{trunc}}s=s$. Moreover, by~\eqref{eq:thb_block} and~\eqref{eq:thb_coarse_coefficients}, $c_F=K_F\hat c_A+K_Tc_T$, that is, $g_F=K_F\hat c_A$, whence $r_{\loc}=\mathbf 0$. Since $\Pi^{\mathrm{trunc}}$ takes values in $\Span\hat\Tau$, it is therefore a projection onto this space.
\end{proof}

As in the HB setting, the residual bounds the error of the step with no unknown constants.

\begin{proposition}[Constant-free $L^\infty$ error bound for THB-splines]\label{prop:thb_linf_bound}
    Let $s\in\Span\Tau$ and $\hat s:=\Pi^{\mathrm{trunc}}s$. Then
    \begin{equation}\label{eq:thb_inf_bound}
        s-\hat s=\sum_{\gamma\in\mathcal C_{\loc}^{\fineonly}}(r_{\loc})_\gamma\,T\gamma
        \qquad\text{and}\qquad
        \|s-\hat s\|_{L^\infty(\Omega)}\le\|r_{\loc}\|_\infty .
    \end{equation}
    Moreover, $s-\hat s$ vanishes identically outside $\bigcup\{\supp T\gamma:\gamma\in\mathcal C_{\loc}^{\fineonly}\}\subseteq\bigcup\{\supp\gamma:\gamma\in\mathcal C_{\loc}^{\fineonly}\}$.
\end{proposition}

\begin{proof}
    Since $r_{\loc}=g_F-K_F\hat c_A$, the identity in~\eqref{eq:thb_inf_bound} is~\eqref{eq:thb_difference} with $\hat v=\hat s$ and $a=\hat c_A$. Each $T\gamma$ is a combination of descendants of $\gamma$, as observed after~\eqref{eq:csub_recursion}, so that $\supp T\gamma\subseteq\supp\gamma$, and the localization follows. Finally, the functions of $\Tau$ are non-negative, see~\eqref{eq:csub_recursion} and the comments after it, and form a partition of unity, so that $\sum_{\gamma\in\mathcal C_{\loc}^{\fineonly}}T\gamma(x)\le1$ for every $x\in\Omega$, whence
    \[
        |s(x)-\hat s(x)|\le\max_{\gamma}|(r_{\loc})_\gamma|\sum_{\gamma\in\mathcal C_{\loc}^{\fineonly}}T\gamma(x)\le\|r_{\loc}\|_\infty .
    \]
\end{proof}

The proof is simpler than that of Proposition~\ref{prop:linf_bound}, which needed the functions carrying the residual to belong to the single level $\ell+1$ in order to use the partition of unity of $\mathcal B_{\ell+1}$; here $\Tau$ is a partition of unity globally. Since Theorem~\ref{thm:adaptive_success} uses only Proposition~\ref{prop:linf_bound}, and the rollback of Section~\ref{S:budget_and_rollback} only that the residual is indexed by $\mathcal C_{\loc}^{\fineonly}$, both remain valid for THB-splines with Proposition~\ref{prop:thb_linf_bound} in place of Proposition~\ref{prop:linf_bound}; the corresponding changes in the algorithms are listed in Section~\ref{S:thb_algorithms}.

\begin{remark}[Two limit cases]\label{rem:thb_rhs}
    If $\mathcal F_{\loc}^{\mathrm{trunc}}=\varnothing$, then $g_F=c_F$ and the least-squares problem coincides with that of the HB setting. If $\mathcal A_{\loc}=\varnothing$, the degenerate configuration of Remark~\ref{rem:no_ancestors}, there is nothing to solve: the step only removes the fine-only functions. In the HB setting, their content is then lost entirely, $r_{\loc}=c_F$. With truncation, instead, the surviving functions of $\mathcal F_{\loc}^{\mathrm{trunc}}$ recover their components along the removed ones, and $r_{\loc}=g_F=c_F-K_Tc_T$. Since $K_T\mathbf 1_{\mathcal F_{\loc}^{\mathrm{trunc}}}=\mathbf 1_{\mathcal C_{\loc}^{\fineonly}}$ by~\eqref{eq:thb_partition_unity_blocks}, this recovery is exact, for instance, when the coefficients of $s$ on $\mathcal C_{\loc}^{\fineonly}\cup\mathcal F_{\loc}^{\mathrm{trunc}}$ are all equal.
\end{remark}

\begin{remark}[Local versus global least squares]\label{rem:thb_locality}
    One could also minimize the whole coefficient residual $\|c-K^{\mathrm{trunc}}\hat c\|_2$ over all $\hat c\in\RR^{\hat n}$, whose entries on $\mathcal A_{\loc}$ we denote by $a$. In the HB setting this changes nothing. By~\eqref{eq:global_refinement_detailed_block}, the coefficients $\hat c_C$ of $\mathcal C_{\loc}^{\mathrm{coarse}}$ appear in a single block of the residual, $c_C-K_Ca-\hat c_C$, and similarly for $\mathcal D_{\loc}$ and $\mathcal H_{\mathrm{rem}}$; these blocks can thus be made zero for every $a$, and the global problem reduces to minimizing $\|c_F-K_Fa\|_2$, as in Theorem~\ref{thm:exact_residual}. With truncation, instead, the coefficients $\hat c_T$ of $\mathcal F_{\loc}^{\mathrm{trunc}}$ appear in two blocks, by~\eqref{eq:thb_block}: $c_F-K_Fa-K_T\hat c_T$ on $\mathcal C_{\loc}^{\fineonly}$ and $c_T-\hat c_T$ on $\mathcal F_{\loc}^{\mathrm{trunc}}$. The minimizer then takes in general $\hat c_T\neq c_T$, reducing the first block at the expense of changing the coefficients of the coarse functions of $\mathcal F_{\loc}^{\mathrm{trunc}}$, and hence the spline on their supports, which extend beyond $\Omega_{\loc}$. The class~\eqref{eq:thb_admissible_class} rules this out by fixing $\hat c_T=c_T$, which leaves the residual $g_F-K_Fa$.
\end{remark}

\begin{remark}[Computational cost]\label{rem:thb_cost}
    Compared with Remark~\ref{rem:complexity}, the least-squares problem is the same, the counter-corrections and the transfer matrix $M_{\ell+1}$ are no longer needed, and the only additional work is the assembly of $K_T$, from $C_\ell^{\ell+1}$ and the matrix $C^{\mathrm{sub}}_\ell$ stored by the hierarchical space, and the product $K_Tc_T$.
\end{remark}

\subsection{Changes in the algorithms}
\label{S:thb_algorithms}

The driver of Algorithm~\ref{alg:coarsen_multilevel} and the adaptive procedure of Algorithm~\ref{alg:adaptive_coarsen} are used unchanged, with $\Tau$ in place of $\HH$. The two local routines require the following modifications, and no other.

\paragraph{In \textnormal{\bf local\_coarsening\_residual} (Algorithm~\ref{alg:estimate_error}).} After the assembly of $K_F$ in Line~4, the following lines are inserted,
\begin{algorithmic}
	\State $E_{\ell+1} \gets C_\ell^{\ell+1}(\mathcal{C}_{\loc}^{\fineonly},\,\cdot\,)\ C^{\mathrm{sub}}_\ell$ \Comment{Un-truncated level-$(\ell+1)$ expansion}
	\State $\mathcal F_{\loc}^{\mathrm{trunc}} \gets \{\beta\in\mathcal{F}^A_{0:\ell} : E_{\ell+1}(\,\cdot\,,\beta)\neq\mathbf 0\}$, \quad $K_T \gets E_{\ell+1}(\,\cdot\,,\mathcal F_{\loc}^{\mathrm{trunc}})$
	\State $g_F \gets c(\mathcal{C}_{\loc}^{\fineonly}) - K_T\, c(\mathcal F_{\loc}^{\mathrm{trunc}})$ \Comment{Modified data, cf.~\eqref{eq:thb_local_vectors}}
\end{algorithmic}
and Lines~5--6 are evaluated at $g_F$ instead of $c(\mathcal C_{\loc}^{\fineonly})$, that is, $\hat c_A\gets K_F^\dagger g_F$ and $r_{\loc}\gets g_F-K_F\hat c_A$. The product in the first line only involves the rows of $C^{\mathrm{sub}}_\ell$ indexed by the parents of the fine-only functions.

\paragraph{In \textnormal{\bf coarsen\_at\_level} (Algorithm~\ref{alg:coarsen_at_level}).} Steps~1--3 are deleted, Step~4, which updates the mesh and the basis, is unchanged, and Step~5 reduces to assembling $\hat c$ by~\eqref{eq:thb_coarse_coefficients}: $\hat c(\mathcal A_{\loc})\gets\hat c_A$, and the coefficients of the surviving indices $\mathcal F_{\loc}^{\mathrm{trunc}}\cup\mathcal F_{\mathrm{rem}}$ are copied from $c$, with no further computation.

Truncation thus simplifies the update at the cost of a slightly larger evaluation: the counter-corrections, and with them the matrices $M_\ell$, disappear from \textbf{coarsen\_at\_level}, while \textbf{local\_coarsening\_residual} gains three lines. During the rollback loop of Algorithm~\ref{alg:adaptive_coarsen} the mesh is not modified, and the candidate set only shrinks, so that the new fine-only functions are a subset of the previous ones; $E_{\ell+1}$ can therefore be computed once, for the initial candidate set, and its rows extracted at each iteration.

\section{Numerical experiments}
\label{sec:numerical_experiments}

The coarsening operators of Sections~\ref{sec:coarsening} and~\ref{sec:thb}, and the adaptive algorithm of Section~\ref{sec:adaptive_coarsening}, are illustrated by two experiments on $\Omega=[0,1]^2$, with cubic splines of maximal smoothness. The first one transfers an adaptive Galerkin solution back along its refinement path, comparing the result with that of the global $L^2$-projection; the second one removes, under a prescribed tolerance, the refinement left behind by a moving front. All the computations were carried out with GeoPDEs~\cite{GV18}.

\subsection{Coarsening along an adaptive path}
\label{S:num_adaptive_path}

We consider the Laplace problem on $\Omega$ with Dirichlet data and source chosen so that the exact solution is the sum of two sharp Gaussians, $u(x)=\sum_{k=1}^{2}\exp(-300|x-a_k|^2)$, with $a_1=(7/12,1/4)$ and $a_2=(3/4,3/4)$. Starting from a uniform mesh of $12\times12$ cells and HB-splines, $N=14$ steps of an adaptive loop driven by a residual-based estimator (cf.~\cite{BuffaGiannelli16}) with the maximum marking strategy, which marks the cells $Q$ whose indicators satisfy $\eta_Q\ge\theta\max_{Q'}\eta_{Q'}$, with $\theta=0.5$, produce a sequence of nested hierarchical meshes, and the Galerkin solutions $s_0,\dots,s_N$ on them, which we call the forward path; the finest mesh, that of $s_N$, has $11223$ degrees of freedom on six levels. The spline $s_N$ is then transferred back along the same path by Algorithm~\ref{alg:coarsen_multilevel}, undoing the marked cells in reverse order, so that each point of a transferred curve results from successive coarsening steps starting at~$s_N$.

We compare the local transfer of Section~\ref{sec:coarsening} with the global $L^2$-projection of $s_N$ onto each coarse space, and with two ablations that reveal the role of its ingredients: the \emph{injection}, which sets $\hat c_A=0$ and applies no counter-correction, so that the coefficients of the removed basis functions are simply discarded; and the transfer \emph{without counter-corrections}, which computes $\hat c_A$ as in Section~\ref{S:local_ls_solution} but leaves the coefficients of $\mathcal C_{\loc}^{\mathrm{coarse}}$ and $\mathcal D_{\loc}$ unchanged. For each of these transfers, $\hat s_k$ denotes the spline it produces on the $k$-th mesh.

\begin{figure}[htb]
	\centering
	\includegraphics[width=\textwidth]{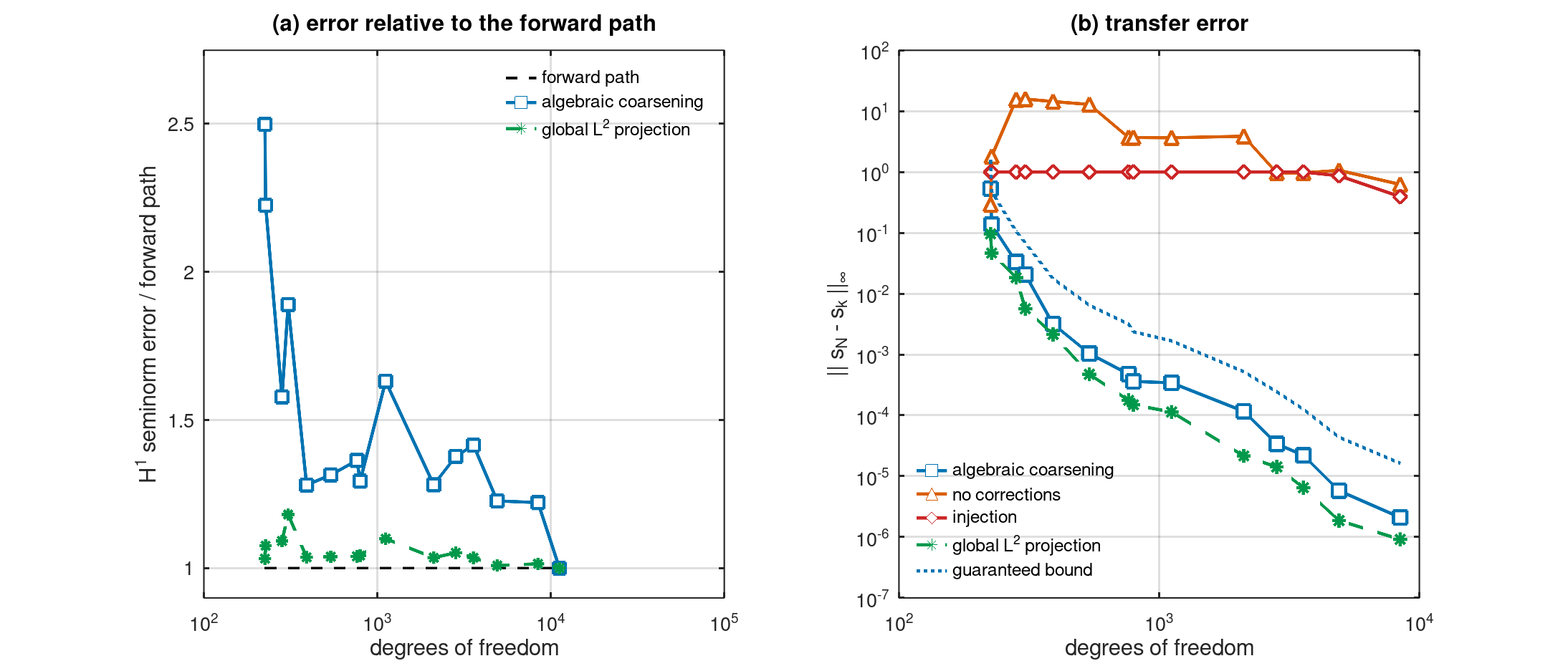}
	\caption{Coarsening along the adaptive path for the Laplace problem whose solution consists of two sharp Gaussians, with bicubic HB-splines of maximal smoothness. Each point of a transferred curve results from successive coarsening steps starting at the finest mesh, which corresponds to the rightmost point. (a) $H^1$-seminorm error of each transfer, divided by that of the forward path. (b) Distance $\|s_N-\hat s_k\|_{L^\infty(\Omega)}$ for each transfer, and bound of the local transfer accumulated along the path.}
	\label{fig:adaptive_path}
\end{figure}

Figure~\ref{fig:adaptive_path}(a) shows the error $|u-\hat s_k|_{H^1(\Omega)}$ of each transfer divided by that of the forward path, $|u-s_k|_{H^1(\Omega)}$. For a single coarsening step, which is the relevant regime when a solve takes place between two coarsening steps, the local transfer multiplies the error of the forward path by $1.22$, whereas the global projection stays within a factor $1.18$ along the whole path. The factor of the local transfer remains moderate over most of the path, so that the transferred splines retain the convergence rate of the forward path, and rises to about $2.5$ only on the two coarsest meshes, after the successive transfers have undone nearly all the refinement. The global projection is more accurate because it may modify every coefficient of the spline; the local transfer gives up part of that accuracy in exchange for leaving the spline untouched away from the coarsened region.

Figure~\ref{fig:adaptive_path}(b) shows the distance $\|s_N-\hat s_k\|_{L^\infty(\Omega)}$ for each transfer, computed exactly, since $\hat s_k$ can be written in the finest basis, together with the bound $e_\infty$ for the local transfer, obtained by adding the bounds of Proposition~\ref{prop:linf_bound} along the path, as Algorithm~\ref{alg:coarsen_multilevel} does across levels. The bound holds at every step, and its effectivity index $e_\infty/\|s_N-\hat s_k\|_{L^\infty(\Omega)}$, that is, the factor by which it overestimates the error, lies between $2.9$ and $7.7$. Both ablations fail by several orders of magnitude, which shows that neither the local least-squares problem nor the counter-corrections can be dispensed with. Finally, the test confirms the locality established in Proposition~\ref{prop:linf_bound}: at a typical step, the local transfer leaves $s_N$ exactly unchanged on $56\%$ of the area of the domain, measured as the set where $\hat s_k$ and $s_N$ coincide up to round-off, whereas there the global projection introduces a perturbation of about $12\%$ of its maximum error.

\subsection{Coarsening behind a moving front}
\label{S:num_moving_front}

In time-dependent simulations with a moving feature, an adaptive mesh remains refined where the feature used to be, and that refinement should be removed. We emulate this situation by moving the two peaks of Section~\ref{S:num_adaptive_path} along a horizontal trajectory. Starting from the same initial mesh, we approximate in turn the functions $u_j:=u(\cdot-d_j)$, $j=1,\dots,4$, with $d_1=(-0.3,0)$, $d_2=(-0.2,0)$, $d_3=(-0.1,0)$ and $d_4=(0,0)$, so that $u_4=u$, without any coarsening in between. For each $j$, five steps of an adaptive loop are performed on the current mesh. Each step computes the $L^2$-projection $u_{j,h}$ of $u_j$ onto the THB-spline space and the indicators $\eta_\beta:=\big(\int_\Omega|u_j-u_{j,h}|^2\beta\big)^{1/2}$, which satisfy $\sum_\beta\eta_\beta^2=\|u_j-u_{j,h}\|^2_{L^2(\Omega)}$ because the THB-splines form a partition of unity; it then marks the basis functions $\beta$ with $\eta_\beta\ge\theta\max_\gamma\eta_\gamma$, where $\theta=0.4$, and refines, for each marked $\beta$, the active cells of the level of $\beta$ in its support, which deactivates $\beta$; see~\cite{GV18}. The spline $s$ to be coarsened is the $L^2$-projection of $u$ onto the THB-spline space of the final mesh, which has $24432$ degrees of freedom on six levels, and its error is $\varepsilon_h:=\|u-s\|_{L^\infty(\Omega)}=2.07\times10^{-7}$.

To measure how much of this refinement is no longer needed, we build a reference mesh by refining for $u$ alone until the same accuracy is reached, with the same loop and $\theta=0.9$, so that each step marks few functions and the loop stops close to that accuracy; it has $15906$ degrees of freedom, also on six levels, and error $0.89\,\varepsilon_h$. We call \emph{excess} the difference between the two meshes, $8526$ degrees of freedom.

We run Algorithm~\ref{alg:adaptive_coarsen} starting from $s$, with THB-splines, for $\tau=\kappa\varepsilon_h$ with $\kappa$ between $0.01$ and $1$ and for the two budget rules of Remark~\ref{rem:budget_distribution}, each tolerance being an independent run. Table~\ref{tab:moving_front} reports the fraction of the excess that is removed, the bound $e_\infty$, and the error $\|s-\hat s\|_{L^\infty(\Omega)}$, whose ratio is the effectivity index of the bound.

\begin{table}[htb]
	\centering
	\caption{Coarsening behind a moving front with THB-splines, one independent run per tolerance $\tau=\kappa\varepsilon_h$ from the same fine state: fraction of the excess removed, bound $e_\infty$ and error $\|s-\hat s\|_{L^\infty(\Omega)}$, both in units of $\varepsilon_h$.}
	\label{tab:moving_front}
	\begin{tabular}{c|ccc|ccc}
		\hline
		& \multicolumn{3}{c|}{greedy} & \multicolumn{3}{c}{equidistributed} \\
				$\kappa$ & excess removed & $e_\infty/\varepsilon_h$ & error$/\varepsilon_h$ & excess removed & $e_\infty/\varepsilon_h$ & error$/\varepsilon_h$ \\
		\hline
	0.01  & 22\% & 0.0100 & 0.00460 & 33\% & 0.00622 & 0.00136 \\
0.03  & 36\% & 0.0300 & 0.00799 & 38\% & 0.0258 & 0.00311 \\
0.1   & 43\% & 0.0999 & 0.0480 & 43\% & 0.0871 & 0.00910 \\
0.3   & 43\% & 0.300 & 0.146 & 49\% & 0.298 & 0.0373 \\
1     & 62\% & 1.00 & 0.395 & 62\% & 0.805 & 0.0952 \\
		\hline
	\end{tabular}
\end{table}

The bound holds in every run. Moreover, the coarsened splines are as accurate as $s$: $\|u-\hat s\|_{L^\infty(\Omega)}$ equals $\varepsilon_h$ to within $0.01\%$ for every tolerance. The coarsening modifies $s$ only along the trail, where its error is much smaller than $\varepsilon_h$, and leaves untouched the neighbourhood of the peaks, where the maximum error is attained. With the greedy rule the whole budget is spent, and the algorithm removes from $22\%$ to $62\%$ of the excess as $\kappa$ grows from $0.01$ to $1$. The compression stalls, however, between $\kappa=0.1$ and $\kappa=0.3$: this rule lets the finest level, where coarsening is more expensive, consume the whole budget.

\begin{figure}[htb]
	\centering
	\includegraphics[width=\textwidth]{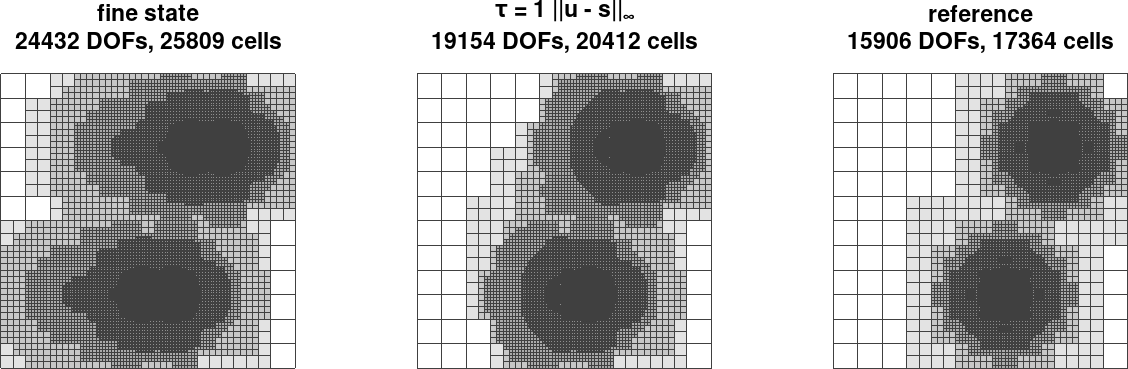}
	\caption{Fine mesh (left), mesh coarsened with tolerance $\tau=\varepsilon_h$ using the greedy budget rule (centre), and reference mesh (right). Darker cells belong to finer levels.}
	\label{fig:moving_front}
\end{figure}

The equidistributed rule reserves part of the budget for the coarser levels. It removes about as much of the excess as the greedy rule or more, $33\%$ instead of $22\%$ for $\kappa=0.01$, with a compression that grows steadily; and for every tolerance both its bound, between $62\%$ and $99\%$ of $\tau$, and its actual error are smaller, the latter by a factor between $2.6$ and $5.3$. This indicates that the coarser levels remove refinement at a lower cost than the finest one. Its bound, however, overestimates its own error by a factor between $4.6$ and $9.6$, instead of $2.1$ and $3.7$: several levels then contribute comparable terms to the sum, whereas the perturbations they introduce need not accumulate in the actual error. A sharper accounting, such as the localized budget of Remark~\ref{rem:local_budget}, would therefore allow even more refinement to be removed. Finally, HB-splines give the same results on the same fine mesh, to all the digits reported in Table~\ref{tab:moving_front}. Figure~\ref{fig:moving_front} compares the fine mesh with the mesh coarsened with $\tau=\varepsilon_h$ using the greedy rule, and with the reference mesh: most of the refinement along the trail has been removed, whereas the refinement around the current position of the peaks is kept.

\section*{Conclusions}
\label{sec:conclusions}

We have presented a purely algebraic coarsening procedure for hierarchical spline spaces, based on a local left inverse of the refinement matrix. A coarsening step is local: it solves a small least-squares problem for the coefficients of the reactivated B-splines and, for HB-splines, applies explicit counter-corrections to the coefficients of the retained B-splines involved in their expansion. It requires neither quadrature nor the solution of a global system. The transfer is exact on the coarse space and leaves the spline unchanged away from the coarsened region. Its residual, available before the mesh is modified, bounds the $L^\infty$-error with no unknown constants, and this bound drives an adaptive multi-level algorithm with a guaranteed tolerance. For THB-splines, truncation makes the counter-corrections unnecessary and modifies only the right-hand side of the local problem.

The experiments show that locality costs little accuracy compared with a global $L^2$-projection, that the bound overestimates the actual error by less than a factor of ten, and that a substantial part of the refinement left behind by a moving front can be removed without affecting the accuracy with respect to the exact solution.

Several directions remain open. The localized budget of Remark~\ref{rem:local_budget} would sharpen the error control when several levels are coarsened in separate regions. Cheaper left inverses of $K_F$, built from quasi-interpolants or from those of tensor-product subdivision matrices (Remark~\ref{rem:alternative_left_inverses}), could further reduce the cost of a step. Natural continuations of this work are bounds in $L^2$ or energy norms, which would involve stability constants, and the analysis of adaptive schemes that alternate refinement and coarsening, in particular for evolution problems.

	\section*{Acknowledgements}
	This work was partially supported by Universidad Nacional del
	Litoral through grant CAI+D-2024 85520240100018LI, and by Consejo Nacional de Investigaciones Científicas y Técnicas through grant PIP 11220250101092CO (Argentina). This support is gratefully acknowledged.

\bibliographystyle{apalike}
\def\cprime{$'$}

\end{document}